\documentclass[a4paper]{amsart}
\usepackage[latin1]{inputenc}
\usepackage{latexsym}
\usepackage{amssymb}
\usepackage{color}
\usepackage{graphics}
\usepackage[all]{xy}

\newcommand{\oprocendsymbol}{\hbox{$\bullet$}}
\newcommand{\oprocend}{\relax\ifmmode\else\unskip\hfill\fi\oprocendsymbol}

\newtheorem{definition}{Definition}[section]
\newtheorem{lemma}[definition]{Lemma}
\newtheorem{theorem}[definition]{Theorem}

\newtheorem{remark}[definition]{Remark}

\begin{document}

\title[Mechanics and Lagrangian submanifolds of presymplectic and Poisson manifolds]
{Time-dependent Mechanics and Lagrangian submanifolds of
presymplectic and Poisson manifolds}

\author[E.\ Guzm\'an, J.C. Marrero]{E.\ Guzm\'an, J.C. Marrero}
\address{E.\ Guzm\'an: Unidad asociada ULL-CSIC Geometr{\'\i}a Diferencial
y Mec\'anica Geom\'etrica, Departamento de Matem\'atica
Fundamental, Facultad de Matem\'aticas, Universidad de la Laguna,
La Laguna, Tenerife, Canary Islands, Spain} \email{eguzman@ull.es}

\thanks{\noindent This work has been partially supported by MEC (Spain) Grants MTM2009-13383,
MTM2009-08166-E and by the Canary Goverment project ACIISI SOLSUBC
200801000238. E. Guzm\'an wishes to thank the CSIC for a
JAE-predoc grant}

\thanks{\noindent {\it Mathematics Subject Classification} (2010): 53D12,
53D17, 70G45, 70H03,  70H05, 70H15.}

\thanks{\noindent {\it Key words and phrases}: Time-dependent Mechanics, presymplectic structures, Poisson structures,
Lagrangian submanifolds, Lagrangian formalism, Hamiltonian
formalism, jet manifolds, Tulczyjew triple}

\begin{abstract}
\noindent A description of time-dependent Mechanics in terms of
Lagrangian submanifolds of
presymplectic and Poisson manifolds is presented. Two new
Tulczyjew triples are discussed. The first one is adapted to the
restricted Hamiltonian formalism and the second one is adapted to
the extended Hamiltonian formalism.
\end{abstract}

\address{J.C. Marrero: Unidad asociada ULL-CSIC Geometr{\'\i}a Diferencial
y Mec\'anica Geom\'etrica, Departamento de Matem\'atica
Fundamental, Facultad de Ma\-te\-m\'a\-ti\-cas, Universidad de la
Laguna, La Laguna, Tenerife, Canary Islands, Spain}
\email{jcmarrer@ull.es}

\maketitle

\tableofcontents


\section{Introduction}
It is well-known that the phase space of velocities of a
mechanical system may be identified with the tangent bundle TQ of
the configuration space Q. Under this identification, the
Lagrangian function is a real $ C^{\infty}$-function  L on TQ and
the Euler-Lagrange equations are
$$\frac{d}{dt}\Big{(}\frac{\partial L}{\partial \dot{q}^{i}}\Big{)} - \frac{\partial L}{\partial q^{i}}
=0, \quad i = 1,..., n= dim Q$$ where $(q^{i}, \dot{q}^{i})$ are
local fibred coordinates on TQ, which represent the positions and
the velocities of the system, respectively.

If the Lagrangian function is hyperregular one may define the
Hamiltonian function $H: T^{*}Q  \longrightarrow \mathbb{R}$ on
the phase space of momenta $T^{*}Q$ and the Euler-Lagrange
equations are equivalent to the Hamilton equations for H
$$\frac{dq^{i}}{dt} = \frac{\partial H}{\partial p_{i}}, \quad \frac{dp_{i}}{dt}
= -  \frac{\partial H}{\partial q^{i}}, \quad i=1, ..., n.$$ Here,
$(q^{i}, p_{i})$ are local fibred coordinates on $T^{*}Q$ which
represent the positions and the momenta of the system,
respectively.

Solutions of the previous Hamilton equations are just the integral
curves of the Hamiltonian vector field $X_{H}$ on $T^{*}Q$ which
is characterized by the condition
$$\iota_{X_{H}}\Omega_{Q} = dH,$$
$\Omega_{Q}$ being the canonical symplectic structure of $T^{*}Q$
(for more details see, for instance, \cite{AbMa,LeRo}).

Lagrangian (Hamiltonian) Mechanics may be also formulated in terms
of Lagrangian submanifolds of symplectic manifolds (see
\cite{Tu1,Tu2}).

In fact, the complete lift $\Omega_{Q}^c$ of $\Omega_{Q}$ to
$T(T^{*}Q)$ defines a symplectic structure on $T(T^{*}Q)$ and, if
on $T^{*}(TQ)$ we consider the canonical symplectic structure
$\Omega_{TQ}$, the canonical Tulczyjew diffeomorphism $A_{Q}:
T(T^{*}Q) \longrightarrow T^{*}(TQ)$ is a symplectic isomorphism.
Moreover, $S_{L}= A_{Q}^{-1}(dL)$ is a Lagrangian submanifold of
the symplectic manifold $(T(T^{*}Q), \Omega_{Q}^{c})$ and the
local equations defining $S_{L}$ as a submanifold of $T(T^{*}Q)$
are just the Euler- Lagrange equations for L.

\begin{figure}[h]
\[
\xymatrix{
&S_{L}\ar[dr]&&S_{H}\ar[dl]&&\\
T^{*}(TQ)\ar[rd]^{\pi_{TQ}}&&T(T^{*}Q)\ar[rr]^{\emph{b}_{\Omega_{Q}}}\ar[ll]_{A_{Q}}\ar[ld]_{T\pi_{Q}}\ar[rd]^{\tau_{T^{*}Q}}&&T^{*}(T^{*}Q)\ar[ld]_{\pi_{T^{*}Q}}\\
&TQ\ar@<1ex>[ul]^{dL}\ar[rr]^{leg_{L}}&&T^{*}Q\ar@<-1ex>[ur]_{dH}&&\\
}
\]
\caption{\it Tulczyjew triple for time-independent Mechanics}
\label{fig1}
\end{figure}
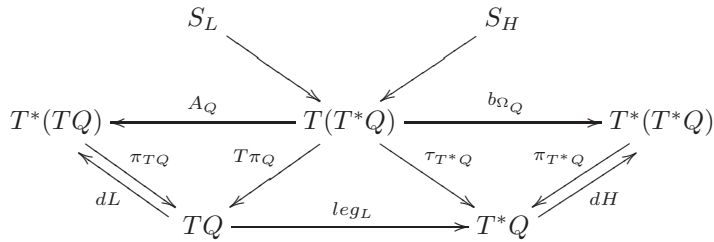

On the other hand, if $H: T^{*}Q  \longrightarrow \mathbb{R}$ is a
Hamiltonian function and $\emph{b}_{\Omega_{Q}}: T(T^{*}Q)
\longrightarrow T^{*}(T^{*}Q)$ is the vector bundle isomorphism
induced by $\Omega_{Q}$ then $\emph{b}_{\Omega_{Q}}$ is an
anti-symplectic isomorphism (when on $T^{*}(T^{*}Q)$ we consider
the canonical symplectic structure $\Omega_{T^{*}Q}$). In
addition, $S_{H}= \emph{b}_{\Omega_{Q}}^{-1}(dH)$ is a Lagrangian
submanifold of $T(T^{*}Q)$ and the local equations defining
$S_{H}$ as a submanifold of $T(T^{*}Q)$ are just the Hamilton
equations for H. Figure \ref{fig1} illustrates the situation

 If the Lagrangian function L is hyperregular then the
Legendre transformation $leg_{L}: TQ\longrightarrow T^{*}Q$ is a
global diffeomorphism and $S_{L} = S_{H}$.

We remark that in the previous construction the following properties hold:
\begin{enumerate}
\item The three spaces $T^{*}(TQ), T(T^{*}Q)$ and $T^{*}(T^{*}Q)$ involved
in the Tulczyjew triple are of the same type, namely, symplectic
manifolds.
\item The two maps $A_{Q}$ and $\emph{b}_{\Omega_{Q}}$ involved in the construction are a symplectic isomorphism and an anti-symplectic
isomorphism, respectively.
\item The Lagrangian and the Hamiltonian functions  are not involved in the definition of the triple. In this sense, the triple is canonical.
\item The dynamical equations (Euler-Lagrange and Hamilton equations) are the local equations defining the Lagrangian submanifolds $S_{L}$ and $S_{H}$ of $T(T^{*}Q)$.
\item The construction may be applied to an arbitrary Lagrangian function (not necessarily regular).
\end{enumerate}
On the other hand, for time-dependent mechanical systems the role
of $TQ$ and $T^{*}Q$ is played by the space of 1-jets $J^{1}\pi$
of local sections of a fibration $\pi: M  \longrightarrow
\mathbb{R}$ (in the Lagrangian formalism) and for the dual bundle
$V^{*}\pi$ to the vertical bundle $V\pi$ to $\pi$ (in the
restricted Hamiltonian  formalism) or for the cotangent bundle
$T^{*}M$ to M (in the extended Hamiltonian formalism). For more
details on these topics, see \cite{LeMaMa,Sa}.

Note that $V^{*}\pi$ is not a symplectic manifold, but a Poisson
manifold.

Several attempts to extend the Tulczyjew triple for time-dependent
mechanical systems have been done. However, although acurrate and
interesting, they all exhibit some defect if we compare with the
original Tulczyjew triple for autonomous mechanical systems. In
fact, in \cite{LeLa} the authors described a Tulczyjew triple for
the particular case when the fibration $\pi: M  \longrightarrow
\mathbb{R}$ is trivial, that is, $M = \mathbb{R}\times Q$ and
$\pi$ is the projection on the first factor. They used the
extended formalism and the spaces involved in the construction
were too big.

Later, in \cite{LeMa}, M. de Le\'on et al discussed a Tulczyjew
triple for the same fibration $pr_{1}:  \mathbb{R}\times Q
\longrightarrow Q$. In this case, the Lagrangian and Hamiltonian
functions are involved in the definition of the triple. In this
construction, they used the notion of the complete lift of a
cosymplectic structure.

On the other hand, in \cite{IgMaPaSo} the authors proposed a
restricted Tulczyjew triple for a general fibration $\pi: M
\longrightarrow \mathbb{R}$. However, the Hamiltonian section is
involved in the construction of the triple.

In this paper, we solve the previous problems and deficiences. In
fact, we will propose two new Tulczyjew triples for time-dependent
mechanical systems. The first one is adapted to the restricted
Hamiltonian formalism and the second one is adapted to the
extended Hamiltonian formalism. In this approach, the role of
symplectic structures in the original Tulczyjew triple is played
by presymplectic  and Poisson structures. Then, symplectic
(anti-symplectic) isomorphisms are replaced by presymplectic and
Poisson (anti-presymplectic and anti-Poisson) isomorphisms. In
addition, Lagrangian submanifolds of symplectic manifolds are
replaced by Lagrangian submanifolds of presymplectic and Poisson
manifolds.

The new Tulczyjew triples follow the same philosophy as the
original one (see sections 4, 5 and compare with properties (1),
(2), (3), (4) and (5) of the original Tulczyjew triple).

We also remark that our second Tulczyjew's triple has some
similarities with the Tulczyjew's triple proposed in
\cite{GrGrUr1} although the spaces involved in the definition of
the triple in \cite{GrGrUr1} are different and the structural
applications between them are not isomorphisms.

The paper is structured as follows. In section 2, we recall some
definitions and results on presymplectic and Poisson structures
which we will be used in the rest of the paper. The Lagrangian and
Hamiltonian formalisms in jet manifolds are discussed in section
3. Sections 4 and 5 contain the results of the paper. In fact, the
restricted and extended Tulczyjew triples for time-dependent
Lagrangian and Hamiltonian systems are presented in sections 4 and
5, respectively. The paper ends with our conclusions and a
description of future research directions.


\section{Presymplectic and Poisson manifolds}

\subsection{Presymplectic manifolds}
$\mbox{}$

In this subsection, we will recall some well-known facts on  presymplectic manifolds.

\begin{definition}\rm
A \emph{presymplectic structure} on a manifold M is a closed 2-form $\omega$ on M. \\
If $\omega$ is a presymplectic structure on M, the couple $(M,
\omega)$ is said to be a  \emph{presymplectic manifold}.
\end{definition}
Moreover, for each $x \in M$, we will denote by $Ker(\omega(x))$
the subspace of the tangent space $T_{x}M$ to M at $x$ given by
$$Ker(\omega(x)) = \{v \in T_{x}M / \iota_{v}\omega(x)=0\}.$$

In other words, $Ker( \omega(x)) = Ker[{ \flat_{\omega}}_{\mid
T_{x}M}],$ where $\flat_{\omega}: TM \longrightarrow T^{*}M$ is
the vector bundle morphism induced by $\omega$.

Note that
$$dim  [Ker (\omega(x))] = dim M - rank (\omega(x)),$$

$rank(\omega(x))$ being the rank of the 2-form $\omega(x)$ which
is an even number.

In the pariticular case when
$$rank(\omega(x)) = dim M, \quad for \quad all  \quad x \in M$$
then the dimension of M is even and the couple $(M, \omega)$ is a
symplectic manifold (see, for instance, \cite{AbMa}).

\begin{definition}\label{lspm}\rm
A submanifold $\emph{C}$ of dimension r  of a presymplectic
manifold $(M, \omega)$ is said to be \emph{ Lagrangian} if
$i^{*}\omega = 0 $ and $$ r =\frac{ rank(\omega(x))}{2} + dim
(T_{x}C \cap Ker (\omega(x))), \quad  for \quad  all  \quad  x \in
C.$$ Here, $i: C \longrightarrow M $ is the canonical inclusion.
\end{definition}
We remark that if $(M, \omega)$ is a symplectic manifold, then one
recovers the classical notion of a Lagrangian submanifold of a
symplectic manifold (see, for instance, \cite{AbMa}).

The notion of a presymplectic map may be introduced in a natural way.

\begin{definition}\rm
A smooth map $\varphi: M  \longrightarrow N$ between two
presymplectic manifolds $(M,\omega_{M})$ and $(N,\omega_{N})$ is
said to be a \emph{presymplectic map} if  $ \varphi^{*}\omega_{N}
= \omega_{M}$.
\end{definition}

Note that if $\flat_{\omega_{M}}: TM \longrightarrow T^{*}M$ and
$\flat_{\omega_{N}}: TN \longrightarrow T^{*}N$ are the bundle
maps induced by $\omega_{M}$ and $\omega_{N}$, respectively, then
$\varphi$ is a presymplectic map if and only if
$${(\flat_{\omega_{M}})}_{\mid T_{x}M}
= (T_{x}\varphi)^{*} \circ { (\flat_{\omega_{N}})}_{\mid T_{\varphi(x)}N} \circ
T_{x}\varphi,$$ for every $x \in M$.

\begin{remark}\rm
A presymplectic structure $\omega$ on a manifold M is a particular
example of a Dirac structure (see \cite{Co}) in such a way that a
Lagrangian submanifold of $(M, \omega)$ is also a Lagrangian
submanifold for the Dirac structure on M  which is induced by the
presymplectic form $\omega$ (see \cite{Uc}). In addition, a
presymplectic map is a backward Dirac map in the sense of Bursztyn
et al (see \cite{BuRa}).
\end{remark}

\subsection{Poisson manifolds}
$\mbox{}$

In this subsection, we will recall some well-known facts on
Poisson manifolds (see, for instance, \cite{GrUr,LiMa}).

\begin{definition}\rm
A \emph{Poisson structure} on a manifold M is a 2-vector $\Lambda$ on M such that
$[\Lambda, \Lambda] = 0 $, where $[.,.]$ is the Schouten-Nijenhuis bracket. \\
If $\Lambda$ is a Poisson structure on M, the couple  $(M,\Lambda)$ is said to be a\emph{ Poisson manifold}.
\end{definition}

A Poisson structure induces a vector bundle morphism
 $\Lambda^{\sharp}: T^{*}M
\longrightarrow TM$ which is given by
$$\Lambda^{\sharp}(\alpha) = \Lambda(\alpha,\ -), \quad for
\quad  \alpha \in T^{*}M.$$

Note that $\Lambda^{\sharp}$ is a skew-symmetric map and, thus,
the dimension of the subspace $\Lambda^{\sharp}(T^{*}_{x}M)$ is
even, for every $x \in M$. Moreover, if  $\Lambda^{\sharp}$ is a
vector bundle isomorphism then the inverse morphism
$(\Lambda^{\sharp})^{-1}: TM\longrightarrow  T^{*}M$ is  just the
vector bundle isomorphism induced by a symplectic structure on M.

\begin{definition}\rm
A  submanifold $\emph{C}$ of a Poisson manifold $(M, \Lambda)$ is said to be
\emph{Lagrangian}  if
$$\Lambda(\alpha, \beta) = 0, \quad  for \quad  all \quad  (\alpha, \beta) \in
(\Lambda^{\sharp})^{-1}(TC)$$ and $$dim (T_{x}C \cap
\Lambda^{\sharp}(T_{x}^{*}M)) = \frac{dim
(\Lambda^{\sharp}(T_{x}^{*}M))}{2}, \quad  for  \quad all \quad x
\in C.$$
\end{definition}

We remark that in the particular case when the map
$\Lambda^{\sharp}: T^{*}M \longrightarrow TM$ is a vector bundle
isomorphism, that is, the Poisson structure is induced by a
symplectic structure on M, then one recovers the classical notion
of a Lagrangian submanifold of a symplectic manifold.
\begin{definition}\rm
A smooth map $\varphi: M  \longrightarrow N$  between two Poisson
manifolds $(M,\Lambda_{M})$ and $(N,\Lambda_{N})$  is said to be a
\emph{Poisson map} if
$$[\Lambda^{2}(T_{x}\varphi)](\Lambda_{M}(x)) = \Lambda_{N}(\varphi(x)),
\quad for \quad each \quad x \in M.$$

Note that $\varphi$ is a Poisson map if and only if
$${(\Lambda^{\sharp}_{N})}_{\mid T_{\varphi(x)}N} =  T_{x}\varphi \circ
{(\Lambda^{\sharp}_{M})}_{\mid T_{x}M} \circ (T_{x}\varphi)^{*},$$
for each $x \in M$.
\end{definition}
\begin{remark}\rm
A Poisson structure $\Lambda$ on a manifold M is a particular
example of a Dirac structure (see \cite{Co}) in such a way that a
Lagrangian submanifold of $(M, \Lambda)$ is also a Lagrangian
submanifold for the Dirac structure on M which is induced by the
Poisson 2-vector $\Lambda$ (see \cite{Uc}). In addition, a Poisson
map is a forward Dirac map in the  sense of Bursztyn et al (see
\cite{BuRa}).
\end{remark}


\section{Lagrangian and Hamiltonian formalisms in jet manifolds}
\setcounter{equation}{0}

In this section, we will recall some definitions and results about
the Lagrangian and Hamiltonian formalisms of Classical Mechanics
in jet manifolds (for more details, see for instance \cite{GrGrUr,
IgMaPaSo, LeMaMa, Sa}).

\subsection{The Lagrangian formalism}
$\mbox{}$

Let $\pi: M \longrightarrow \mathbb{R}$ be a fibration, where M is
a manifold of dimension n+1.

Denote by $J^{1}\pi$  the (2n+1)-dimensional manifold of 1-jets of
local sections of $\pi$. $J^{1}\pi$ is an affine bundle modelled
over the vertical bundle  $V\pi$ of $\pi$.  It can be shown that
exits a canonical  identification between $J^{1}\pi$ and the
subset of $TM$ given by $\{v \in TM / \eta(v)=1\}$, where
\emph{$\eta=\pi^{*}(dt)$}. Thus, $J^{1}\pi$  is an embedded
submanifold of TM. In the same way, $V\pi$  is the vector
subbundle of TM given by $ \{v \in TM / \eta(v)=0\}$.

If $(t, q^{i})$ are local coordinates on M which are adapted to
the fibration $\pi$, then we can consider the corresponding local
coordinates   $(t, q^{i}, \dot{q}^{i})$ on  $J^{1}\pi$ and $V\pi$.

We will denote by $\pi_{1,0}: J^{1}\pi  \longrightarrow M$ and
$\pi_{1}: J^{1}\pi \longrightarrow \mathbb{R}$ the canonical
projections and by $\eta_{1}$ the 1-form on $J^{1}\pi$ given by
$\eta_{1}=(\pi_{1})^{*}(dt)$.

Given the fibration $\pi$, a \emph{Lagrangian function} is a
function $L \in C^{\infty}(J^{1}\pi)$, that is, $L:J^{1}\pi
\longrightarrow \mathbb{R}$.

Given two points $x, y \in M$ we define the manifold of infinite
piecewise differentiable local sections which connect x and y as
$$C^{\infty}(x,y) = \{ c: [0,1] \longrightarrow M / c \mbox{ is a local section of } \pi,
c(0)=x \mbox{ and } c(1)=y \}.$$

We define the functional  $\mathcal{J}: C^{\infty} (x,y) \longrightarrow \mathbb{R}$ by
$$c \rightarrowtail \mathcal{J}(c)= \int^{1}_{0}L(j^{1}c(t))dt. $$

Here, $j^{1}c: [0,1]  \longrightarrow J^{1}\pi$ is the jet
prolongation of the curve c.

The Hamilton principle states that a curve $c \in C^{\infty}
(x,y)$ is a motion of the Lagrangian system defined by L if and
only if c is a critical point on $\mathcal{J}$, i.e.,
$d\mathcal{J}(c)(X)=0$ for all $X\in T_{c} C^{\infty} (x,y)$ which
is equivalent to the condition
\begin{equation}\label{a1}
\frac{d}{dt}(\frac{\partial L}{\partial \dot{q}^{i}})- \frac{\partial L}{\partial q^{i}} =0, \forall i .
\end{equation}

In other words, c satisfies the Euler-Lagrange equations.

\subsection{The Hamiltonian formalism}\label{Ham-formalism}
$\mbox{}$

Denote by $V^{*}\pi$ the dual bundle to the vertical bundle to
$\pi$ and by $\mu: T^{*}M  \longrightarrow V^{*}\pi$ the canonical
projection. We have that $T^{*}M$ is an affine bundle over
$V^{*}\pi$ of rank 1 modelled over the trivial vector bundle
$pr_{1}:  V^{*}\pi \times \mathbb{R} \longrightarrow V^{*}\pi$ (an
AV-bundle in the terminology of \cite{GrGrUr}).

In this setting, a  \emph{Hamiltonian section} is a section
$h:V^{*}\pi \longrightarrow T^{*}M$ of $\mu: T^{*}M
\longrightarrow V^{*}\pi$.

If $(t,q^{i},p,p_{i})$ (respectively, $(t,q^{i},p_{i})$) are local
coordinates on $T^{*}M$ (respectively, $ V^{*}\pi$) we have that
$$\mu (t,q^{i},p,p_{i})= (t,q^{i},p_{i}),\quad h (t,q^{i},p_{i})=(t,q^{i}, -H(t,q^{i},p_{i}),p_{i}).$$

Denote by $\Omega_{M}$ the canonical symplectic structure of
$T^{*}M$. Then, we can obtain  a cosymplectic structure
$(\Omega_{h}, \eta_{1}^{*}) $  on $V^{*}\pi$, where
$$\Omega_{h}= h^{*}\Omega_{M}\in \Omega^{2}(V^{*}\pi), \quad
\eta_{1}^{*} = (\pi^{*}_{1})^{*}(dt)\in \Omega^{1}(V^{*}\pi).$$

Here, $\pi_{1}^{*}: V^{*}\pi \longrightarrow \mathbb{R}$ is the
canonical projection. Note that
$$\Omega_{h}=dq^{i}\land dp_{i}
+dH\land dt, \quad \eta_{1}^{*}= dt.$$
 Thus, we can construct the
Reeb vector field of $(\Omega_{h}, \eta_{1}^{*}) $, which is
characterized by the following conditions
$$\iota_{R_{h}}\Omega_{h} =0,\quad \iota _{R_{h}}\eta_{1}^{*} =1.$$

The local expression of $R_{h}$ is
\begin{equation}\label{Rh}
R_{h} =  \frac{\partial }{\partial t} +  \frac{\partial H}{\partial p_{i}}\frac{\partial }{\partial q^{i}}
-  \frac{\partial H}{\partial q^{i}}\frac{\partial }{\partial p_{i}}
\end{equation}
and, therefore, the integral curves of $R_{h}$  are the solutions
of the \emph{ Hamilton  equations}:
\begin{equation}\label{a2}
\frac{dq^{i}}{dt}= \frac{\partial H}{\partial p_{i}},\quad
\frac{dp_{i}}{dt}= -\frac{\partial H}{\partial q^{i}}, \quad
\forall i.
\end{equation}
This is \emph{the restricted formalism} for time-dependent Hamiltonian Mechanics.\\

Next, we will present \emph{the extended formalism}.

The AV-bundle $\mu: T^{*}M \longrightarrow V^{*}\pi$ is a
principal $ \mathbb{R}$-bundle. We will denote by $V_{\mu} \in
\frak{X}(T^{*}M) $ the infinitesimal generator of the action of $
\mathbb{R}$ on $T^{*}M$. Then, there exists
 a one-to-one correspondence between the space $\Gamma(\mu)$ of
 sections of $\mu$ and the set $\{F_{h} \in C^{\infty}(T^{*}M ) /
 V_{\mu}(F_{h})=1\}$.
Thus, the Hamiltonian section $h: V^{*}\pi \longrightarrow T^{*}M$
induces a real function $F_{h} \in C^{\infty}(T^{*}M )$ such that
$V_{\mu}(F_{h})=1$. The local expression of $F_{h}$ is
\begin{equation}\label{Fh}
F_{h}(t,q^{i},p,p_{i}) =  p + H(t,q^{i},p_{i}).
\end{equation}
Note that $V_{\mu} = \displaystyle {\frac{\partial }{\partial
p}}$.
\begin{remark}\label{princ-connec}{\rm
We remark that $dF_{h}$ is invariant under the action of
$\mathbb{R}$ on $T^*M$ and, thus, it defines a connection $1$-form
on the principal $\mathbb{R}$-bundle $\mu: T^*M \to V^*\pi$. }
\hfill{$\Diamond$}
\end{remark}
 Now, we can consider the Hamiltonian vector field
$\mathcal{H}_{F_{h}}^{\Omega_{M}}$ of $F_{h}$ with respect to the
canonical symplectic structure $\Omega_{M}$.  The local expression
of  $\mathcal{H}_{F_{h}}^{ \Omega_{M}}$ is
\begin{equation}\label{hf}
\mathcal{H}_{F_{h}}^{\Omega_{M}}= \frac{\partial}{\partial t} -\frac{\partial H}{\partial t}\frac{\partial}{\partial p} + \frac{\partial H}{\partial p_{i}}\frac{\partial}{\partial q^{i}}  - \frac{\partial H}{\partial
q^{i}}\frac{\partial }{\partial p_{i}} .
\end{equation}
So, it is clear that $\mathcal{H}_{F_{h}}^{\Omega_{M}}$ is
$\mu$-projectable over $R_{h}$.

In addition, the integral curves of $\mathcal{H}_{F_{h}}^{\Omega_{M}}$ satisfy the following equations
\begin{equation}\label{ec1}
\frac{dq^{i}}{dt}=\frac{\partial H}{\partial p_{i}},  \quad  \frac{dp_{i}}{dt}= -\frac{\partial H}{\partial q^{i}},\quad       i \in \{1, \cdots, m\}
\end{equation}
\noindent and, moreover,
\begin{equation}\label{ec2}
\frac{dp}{dt}= -\frac{\partial H}{\partial t}
\end{equation}
 (\ref{ec1}) are the Hamilton equations and using  (\ref{ec2}) we deduce
 that in time-dependent Mechanics the Hamiltonian energy is not, in general,
 a constant of the motion (for more details, see the following
 subsection 3.3).


\subsection{The equivalence between the Lagrangian and Hamiltonian formalisms}
$\mbox{}$

We are going to introduce the Legendre transformations  for the
restricted and extendend formalisms.

The \emph{extended Legendre transformation} $ Leg_{L}: J^{1}\pi
\longrightarrow T^{*}M$  is given by $
(Leg_{L})(v)(X) = L(v)\eta(X) +{ \frac{d}{dt}}_{\mid t=0}L(v+t(X-\eta(X)v)),$ for $v\in J^{1}\pi$ and $X \in T_{x}M$, with $x = \pi_{1,0}(v).$

The \emph{restricted Legendre transformation} $leg_{L}: J^{1}\pi
\longrightarrow V^{*}\pi$ is defined by  $leg_{L}= \mu \circ
Leg_{L}.$

The local expression of these transformations is
\begin{equation}\label{ec3}
Leg_{L}(t,q^{i},\dot{q}^{i})= (t,q^{i}, L-\dot{q}^{i}\frac{\partial L}{\partial\dot{q}^{i}},\frac{\partial L}{\partial\dot{q}^{i}}), \quad leg_{L}(t,q^{i},\dot{q}^{i})= (t,q^{i},\frac{ \partial L}{\partial\dot{q}^{i}}).
\end{equation}

\vspace{0.1cm} The Lagrangian function L is said to be \emph{ regular} if and only if  for each canonical coordinate system
 $(t, q^{i}, \dot{q}^{i})$  in $J^{1}\pi$, the Hessian matrix
 $ W_{ij} = (\frac{\partial^{2}L}{\partial \dot{q}^{i}\partial \dot{q}^{j}})$ is non-singular.

From (\ref{ec3}), we deduce that  the following statements are equivalent:
\begin{itemize}
\item L is regular.
\item $leg_{L}: J^{1}\pi \longrightarrow V^{*}\pi$  is a local diffeomorphism.
\item $ Leg_{L}: J^{1}\pi \longrightarrow T^{*}M$ is an immersion.
\end{itemize}

\vspace{0.1cm} The Lagrangian function L is said to be
\emph{hyperregular}  if the restricted  Legendre transformation is
a global diffeomorphism. Then, we  obtain a Hamiltonian section $h
= Leg_{L}\circ leg_{L}^{-1}$.  Moreover, if we consider the vector
field $R_L$ on $J^{1}\pi$ given by
\[
R_L(v) = (T_{leg_L(v)}leg_{L}^{-1})(R_h(leg_L(v))), \; \; \mbox{
for } v\in J^1\pi,
\]
then $R_L$ is a second order differential equation on $J^1\pi$ and
the trajectories of $R_L$ are just the solutions of the
Euler-Lagrange equations for $L$. $R_L$ is called \emph{the
Euler-Lagrange vector field} for $L$ and its local expression is
\begin{equation}\label{R-L}
R_L = \displaystyle \frac{\partial}{\partial t} + \dot{q}^i
\frac{\partial}{\partial q^i} + W^{ij}(\frac{\partial L}{\partial
q^i} - \dot{q}^k \frac{\partial^2 L}{\partial \dot{q}^i
\partial q^k} - \frac{\partial^2 L}{\partial t \partial
\dot{q}^i})\frac{\partial}{\partial \dot{q}^j},
\end{equation}
where $(W^{ij})$ is the inverse matrix of $(W_{ij}) =
(\frac{\partial^{2}L}{\partial \dot{q}^{i}\partial \dot{q}^{j}})$.

Using the above facts, we deduce that if $\sigma:
\mathbb{R}\longrightarrow M$ is a solution of the Euler-Lagrange
equations for L then $leg_{L} \circ j^{1}\sigma:
\mathbb{R}\longrightarrow V^{*}\pi$ is a solution of the Hamilton
equations for h and, conversely, if $\tau:
\mathbb{R}\longrightarrow V^{*}\pi$ is a solution of the Hamilton
equations for h then $leg_{L}^{-1} \circ \tau:
\mathbb{R}\longrightarrow J^{1}\pi$ is a prolongation of a
solution $\sigma$ of the Euler-Lagrange equations for L.


\section{Restricted Tulczyjew's  triple}

\subsection{The Lagrangian formalism}
\setcounter{equation}{0}
$\mbox{}$

Let N be a smooth manifold. We will denote by  $A_{N}: T(T^{*}N)
\longrightarrow T^{*}(TN)$ the canonical Tulczyjew diffeomorphism
associated with the manifold N which is given locally by (see
\cite{Tu2})
$$A_{N}(q^{i}, p_{i}; \dot{q}^{i}, \dot{p}_{i})= (q^{i}, \dot{q}^{i}; \dot{p}_{i}, p_{i}).$$
Here $(q^{i})$ are local coordinates on N and $(q^{i}, p_{i})$
(respectively, $(q^{i}, p_{i}; \dot{q}^{i}, \dot{p}_{i}))$ are the
corresponding local coordinates on $T^{*}N$ (respectively,
$T(T^{*}N)$.

Now, suppose that $\pi: M \longrightarrow \mathbb{R}$ is a fibration. Then, we may define a smooth map
$$\psi: T^{*}(J^{1}\pi)  \longrightarrow T(V^{*}\pi)$$
as follows. Let $\alpha_{v}$ be a 1-form at the point $v \in J^{1}\pi \subseteq TM$. Then,
$$ \psi (\alpha_{v}) = T\mu (A_{M}^{-1}(\tilde{\alpha}_{v})),$$
with $\tilde{\alpha}_{v} \in T_{v}^{*}(TM)$ such that
$\tilde{\alpha}_{v\vert T_{v}(J^{1}\pi)}= \alpha_{v}$ and $\mu:
T^{*}M  \longrightarrow V^{*}\pi$ being the canonical projection.

$\psi$ is well-defined. In fact, the local expression of $\psi$ is
\begin{equation}\label{ec4}
\psi (t,q^{i}, \dot{q}^{i};p_{t}, p_{q^{i}}, p_{\dot{q}^{i}})= (t,q^{i},p_{\dot{q}^{i}}; 1,  \dot{q}^{i}, p_{q^{i}}).
\end{equation}
In particular, $\psi $ take values in the submanifold
$J^{1}\pi^{*}_{1}$ of $T(V^{*}\pi).$ Thus, we may consider the map
$$ \psi: T^{*}(J^{1}\pi)  \longrightarrow J^{1}\pi^{*}_{1}.$$

It is clear that $\psi$ is not a diffeomorphism  (see
(\ref{ec4})). In order to obtain a diffeomorphism,  we consider
the vector subbundle  $ \langle\eta_{1}\rangle$ over $J^{1}\pi$ of
$T^{*}(J^{1}\pi)$ with rank 1 which is generated by the 1-form
$\eta_{1}$ and the quotient vector bundle $ T^{*}(J^{1}\pi)/
\langle\eta_{1}\rangle$  over $J^{1}\pi$. Local coordinates on $
T^{*}(J^{1}\pi)/ \langle\eta_{1}\rangle$ are $(t, q^{i},
\dot{q}^{i}; p_{q^{i}}, p_{\dot{q}^{i}})$. In  addition, it is
easy to prove that there exits a diffeomorphism $\tilde{\psi}:
T^{*}(J^{1}\pi)/ \langle\eta_{1}\rangle \longrightarrow
J^{1}\pi^{*}_{1}$ such that the following diagram is commutative
\[
\xymatrix{
T^{*}(J^{1}\pi)\ar[d]^{\pi_{T^{*}(J^{1}\pi)}}\ar[rr]^{\psi}&&J^{1}\pi^{*}_{1}\\
T^{*}(J^{1}\pi)/ \langle\eta_{1}\rangle\ar[urr]_{\tilde{\psi}}&&
 }
\]
where $\pi_{T^{*}(J^{1}\pi)}$ is the canonical projection. In
fact, the local expression of $\tilde{\psi}$ is
$$\tilde{\psi}(t,q^{i}, \dot{q}^{i}; p_{q^{i}}, p_{\dot{q}^{i}}) = (t,q^{i},  p_{\dot{q}^{i}};\dot{q}^{i}, p_{q^{i}}).$$

We will denote by $A_{\pi}:  J^{1}\pi^{*}_{1}\longrightarrow
T^{*}(J^{1}\pi)/ \langle\eta_{1}\rangle$ the inverse of
$\tilde{\psi}$. $A_{\pi}$ will be called \emph{ the canonical
Tulczyjew diffeomorphsim associated with the fibration $\pi$}. The
local expression of $A_{\pi}$ is
\begin{equation}\label{ec5}
A_{\pi} (t, q^{i}, p_{i}; \dot{q}^{i}, \dot{p}_{i}) = (t, q^{i},\dot{q}^{i}; \dot{p}_{i},  p_{i}).
\end{equation}

Let $\Omega_{J^{1}\pi}$ be the canonical symplectic structure of
$T^{*}(J^{1}\pi)$ and $\Lambda_{J^{1}\pi}$ be the corresponding
Poisson structure.

In local coordinates $(t, q^{i}, \dot{q}^{i};p_{t}, p_{q^{i}},
p_{\dot{q}^{i}})$ on $T^{*}(J^{1}\pi)$, we have that
$$\Omega_{J^{1}\pi}= dt\land dp_{t}+ dq^{i}\land dp_{q^{i}} + d\dot{q}^{i} \land dp_{\dot{q}^{i}},$$
$$\Lambda_{J^{1}\pi}=  \frac{\partial }{\partial t}\land \frac{\partial }{\partial p_{t}}+  \frac{\partial }{\partial q^{i}}\land \frac{\partial }{\partial p_{q^{i}}}+  \frac{\partial }{\partial \dot{q}^{i}}\land \frac{\partial }{\partial p_{\dot{q}^{i}}}.$$

On the other hand, the vertical bundle of the canonical projection
$\pi_{T^{*}(J^{1}\pi)}:  T^{*}(J^{1}\pi)\longrightarrow
T^{*}(J^{1}\pi)/ \langle\eta_{1}\rangle$ is generated by the
vertical lift $\eta_{1}^{v}$ of the 1-form $\eta_{1}$ on
$J^{1}\pi$. Note that
$$\eta_{1}^{v}=  \frac{\partial }{\partial p_{t}}.$$
Thus, it is clear that
$$ \displaystyle{\mathcal{L}_{\eta_{1}^{v}} \Lambda_{J^{1}\pi}=0}$$
and, therefore, $\Lambda_{J^{1}\pi}$ is $
\pi_{T^{*}(J^{1}\pi)}$-projectable over a Poisson structure
$\widetilde{\Lambda}_{J^{1}\pi}$ on $T^{*}(J^{1}\pi)/
\langle\eta_{1}\rangle$. In fact,
\begin{equation}\label{ec6}
\widetilde{\Lambda}_{J^{1}\pi }=  \frac{\partial }{\partial q^{i}}\land \frac{\partial }{\partial p_{q^{i}}}+  \frac{\partial }{\partial \dot{q}^{i}}\land \frac{\partial }{\partial p_{\dot{q}^{i}}}.
\end{equation}
The corank of the Poisson structure $\widetilde{\Lambda}_{J^{1}\pi
}$ is 1.

Now, consider the canonical Poisson structure $\Lambda_{V^{*}\pi}$
on $V^{*}\pi$. $\Lambda_{V^{*}\pi}$ is characterized by the
following conditions
$$\Lambda_{V^{*}\pi}(d\widehat{X},d\widehat{Y}) = - \widehat{[X, Y]}, \quad \Lambda_{V^{*}\pi}(d(f \circ \pi_{1,0}^{*}),d\widehat{Y})= Y(f) \circ \pi_{1,0}^{*}, \quad \Lambda_{V^{*}\pi}(d(f \circ \pi_{1,0}^{*}),d(g \circ \pi_{1,0}^{*}))= 0 $$
for $X, Y$ $ \pi$-vertical vector fields on M and f, g $\in
C^{\infty}(M)$, where $\pi_{1,0}^{*}: V^{*}\pi \longrightarrow M$
is the canonical projection. Here, $\widehat{Z}$ is the linear
function on $V^{*}\pi$ which is induced by a $\pi$-vertical vector
field Z on M, that is,
$$\widehat{Z}(\alpha)=\alpha(Z(\pi_{1,0}^{*}(\alpha))),\quad \forall \alpha \in V^{*}\pi.$$
If $(t, q^{i}, p_{i})$ are local coordinates on $V^{*}\pi$ then
$$\Lambda_{V^{*}\pi} =  \frac{\partial }{\partial q^{i}}\land \frac{\partial }{\partial p_{i}}.$$
Next, let $\Lambda_{V^{*}\pi}^{c}$ be the complete lift of
$\Lambda_{V^{*}\pi}$ to $T(V^{*}\pi)$. $\Lambda_{V^{*}\pi}^{c}$ is
a Poisson structure on $T(V^{*}\pi)$. Note that  the local
expression of $\Lambda_{V^{*}\pi}^{c}$ is
$$\Lambda_{V^{*}\pi}^{c} =   \frac{\partial }{\partial q^{i}}\land \frac{\partial }{\partial \dot{p}_{i}}+ \frac{\partial }{\partial \dot{q}^{i}}\land \frac{\partial }{\partial p_{i}} .$$
On the other hand, $J^{1}\pi_{1}^{*}$ is a embedded submanifold of
$T(V^{*}\pi)$. In fact, if $(t, q^{i}, p_{i}; \dot{t},
\dot{q}^{i}, \dot{p}_{i})$ are local coordinates on $T(V^{*}\pi)$
then the local equation definning $J^{1}\pi_{1}^{*}$ as a
submanifold of $T(V^{*}\pi)$ is $\dot{t}=1$.

Thus, the restriction $\Lambda_{J^{1}\pi_{1}^{*}}$ to
$J^{1}\pi_{1}^{*}$ of $\Lambda_{V^{*}\pi}^{c}$ is tangent to
$J^{1}\pi_{1}^{*}$ and, furthemore, $\Lambda_{J^{1}\pi_{1}^{*}}$
defines a Poisson structure on $J^{1}\pi_{1}^{*}$.

If $(t, q^{i}, p_{i}, \dot{q}^{i}, \dot{p}_{i})$ are local coordinates on $J^{1}\pi_{1}^{*}$, we have that
\begin{equation}\label{ec7}
 \Lambda_{J^{1}\pi_{1}^{*}}= \frac{\partial }{\partial q^{i}}\land \frac{\partial }{\partial \dot{p}_{i}}+ \frac{\partial }{\partial \dot{q}^{i}}\land \frac{\partial }{\partial p_{i}}.
\end{equation}
Therefore, $ \Lambda_{J^{1}\pi_{1}^{*}}$ is a Poisson structure of
corank 1.

In addition, from  (\ref{ec5}),  (\ref{ec6}) and  (\ref{ec7}), we deduce
\begin{theorem}
$A_{\pi}$ is a Poisson isomorphism between the Poisson manifolds
$(J^{1}\pi_{1}^{*} ,\Lambda_{J^{1}\pi_{1}^{*}})$ and \\ $(
T^{*}(J^{1}\pi)/ \langle\eta_{1}\rangle,
\widetilde{\Lambda}_{J^{1}\pi })$.
\end{theorem}
The space $\frac{ T^{*}(J^{1}\pi)}{ \langle\eta_{1}\rangle}$ is a
vector bundle over $J^{1}\pi$ with vector bundle projection
$\tilde{\pi}_{J^{1}\pi}: \frac{ T^{*}(J^{1}\pi)}{
\langle\eta_{1}\rangle }\longrightarrow J^{1}\pi$. Moreover, we
can consider the jet prolongation $j^{1}\pi_{1,0}^{*}:
J^{1}\pi_{1}^{*}  \longrightarrow J^{1}\pi$ of the bundle map
$\pi_{1,0}^{*}:V^{*}\pi  \longrightarrow M $. We have that
$$\tilde{\pi}_{J^{1}\pi}(t, q^{i},\dot{q}^{i}; p_{q^{i}}, p_{\dot{q}^{i}}) = (t, q^{i},\dot{q}^{i}).$$
Therefore, it is clear that $\tilde{\pi}_{J^{1}\pi} \circ A_{\pi}
= j^{1}\pi_{1,0}^{*}$.

On the other hand, as we know, $J^{1}\pi_{1}^{*}$ is an affine
bundle over $V^{*}\pi$ which is modelled over the vertical bundle
to $\pi_{1}^{*}: V^{*}\pi  \longrightarrow \mathbb{R}$. We will
denote by $(\pi_{1}^{*})_{1,0}: J^{1}\pi_{1}^{*} \longrightarrow
V^{*}\pi$ the affine bundle projection. It follows that
$(\pi_{1}^{*})_{1,0}(t, q^{i},p_{i}; \dot{q}^{i}, \dot{p}_{i}) =
(t, q^{i}, p_{i})$.

The following commutative diagram illustrates the above situation
\[
\xymatrix{ T^{*}(J^{1}\pi)/
\langle\eta_{1}\rangle\ar[rd]_{\tilde{\pi}_{J^{1}\pi}}&&J^{1}\pi^{*}_{1}\ar[ll]_{A_{\pi}}\ar[rd]^{(\pi_{1}^{*})_{1,0}}\ar[ld]^{j^{1}(\pi_{1,0}^{*})}&&\\
& J^{1}\pi&& V^{*}\pi&&
 }
\]
Now, suppose that $L: J^{1}\pi  \longrightarrow \mathbb{R}$  is a
Lagrangian function. Then, the differential of L induces a section
of the vector bundle $\tilde{\pi}_{J^{1}\pi}: T^{*}(J^{1}\pi)/
\langle\eta_{1}\rangle \longrightarrow J^{1}\pi$ which we will
denote by
$$\widetilde{dL}: J^{1}\pi  \longrightarrow T^{*}(J^{1}\pi)/ \langle\eta_{1}\rangle.$$
We have that
\begin{equation}\label{dl}
\widetilde{dL}(t, q^{i}, \dot{q}^{i}) = ( t, q^{i}, \dot{q}^{i};  \frac{\partial L }{\partial q^{i}},  \frac{\partial L }{\partial \dot{q}^{i}}).
\end{equation}
Furthemore, it is easy to prove that $\widetilde{dL}(J^{1}\pi)$ is
a Lagrangian submanifold of the Poisson manifold $(
T^{*}(J^{1}\pi)/ \langle\eta_{1}\rangle,
\widetilde{\Lambda}_{J^{1}\pi })$. In fact,
$$(\tilde{\Lambda}_{J^{1}\pi}^{\sharp})^{-1} (T(\widetilde{dL}(J^{1}\pi))) = \Big{\langle}\Big{ \{}dp_{q^{j}} - \frac{\partial^{2}L}{\partial q^{i}\partial q^{j}}dq^{i} - \frac{\partial^{2}L}{\partial \dot{q}^{i}\partial q^{j}}d\dot{q}^{i}, \quad
dp_{\dot{q}^{k}} - \frac{\partial^{2}L}{\partial
\dot{q}^{k}\partial q^{l}}dq^{l} - \frac{\partial^{2}L}{\partial
\dot{q}^{k}\partial \dot{q}^{l}}d\dot{q}^{l}\Big{\}}\Big{\rangle}
$$
\noindent and
\begin{eqnarray*}
  T(\widetilde{dL}(J^{1}\pi)) \cap \tilde{\Lambda}_{J^{1}\pi}^{\sharp}\Big{ (}T^{*}\Big{(}\frac{T^{*}(J^{1}\pi)}{\langle \eta_{1} \rangle}\Big{)}\Big{)} &=&\Big{\langle}\Big{ \{ } \frac{\partial}{\partial q^{i}}+ \frac{\partial^{2}L}{\partial q^{i}\partial \dot{q}^{j}}\frac{\partial}{\partial p_{q^{i}}}+  \frac{\partial^{2}L}{\partial \dot{q}^{i}\partial q^{j}}\frac{\partial}{\partial p_{\dot{q}^{i}}},\\
  & & \quad  \frac{\partial}{ \partial \dot{ q}^{k} }+ \frac{\partial^{2}L}{\partial \dot{q}^{k}\partial q^{l}}\frac{\partial}{\partial p_{q^{l}}}+  \frac{\partial^{2}L}{\partial \dot{q}^{k}\partial \dot{q}^{l}}\frac{\partial}{\partial p_{\dot{q}^{l}}} \Big{\}}\Big{\rangle}
\end{eqnarray*}
\noindent which implies that
$$\widetilde{\Lambda}_{J^{1}\pi}(\alpha, \beta) = 0, \quad \forall \alpha, \beta \in (\widetilde{\Lambda}_{J^{1}\pi}^{\sharp})^{-1}(T(\widetilde{dL}(J^{1}\pi))),$$
$$dim \Big{ (T}_{\widetilde{dL}(z)}(\widetilde{dL}(J^{1}\pi)) \cap \widetilde{\Lambda}_{J^{1}\pi}^{\sharp}\Big{(}T^{*}_{\widetilde{dL}(z)}\Big{(}\frac{T^{*}(J^{1}\pi)}{\langle \eta_{1} \rangle}\Big{)}\Big{)}\Big{)} = \frac{dim \Big{(}\widetilde{\Lambda}_{J^{1}\pi}^{\sharp}\Big{(}T_{\widetilde{dL}(z)}^{*}\Big{(}\frac{T^{*}(J^{1}\pi)}{\langle \eta_{1} \rangle} \Big{)} \Big{)}  \Big{)}  }{2} =2n, $$
\noindent $\forall z \in J^{1}\pi.$

Thus, since $A_{\pi}$ is a Poisson isomorphism, we deduce that
$S_{L} = A_{\pi}^{-1}(\widetilde{dL}(J^{1}\pi))$ is a
\emph{Lagrangian submanifold} of the Poisson manifold
$(J^{1}\pi_{1}^{*} ,\Lambda_{J^{1}\pi_{1}^{*}})$.

On the other hand, we will denote by $ leg_{L}: J^{1}\pi
\longrightarrow V^{*}\pi$ the restricted Legendre transformation
associated with L. Then, we have the following result.
\begin{theorem}
\begin{enumerate}
\item Let $\sigma: \mathbb{R} \longrightarrow M $ be a local section of $\pi$. $\sigma$ is a solution of the Euler-Lagrange equations for L if and only if
$$A_{\pi}^{-1}\circ \widetilde{dL} \circ j^{1}\sigma = j^{1}(leg_{L} \circ j^{1}\sigma).$$
\item The local equations which define to $S_{L}$ as a Lagrangian submanifold of the Poisson ma\-nifold $(J^{1}\pi_{1}^{*} ,\Lambda_{J^{1}\pi_{1}^{*}})$ are just the Euler-Lagrange equations for L.
\end{enumerate}
\end{theorem}

\underline{Proof}: A local computation, using (\ref{a1}),
(\ref{ec3}) and (\ref{ec5})  proves the result. \hfill{$\square$}

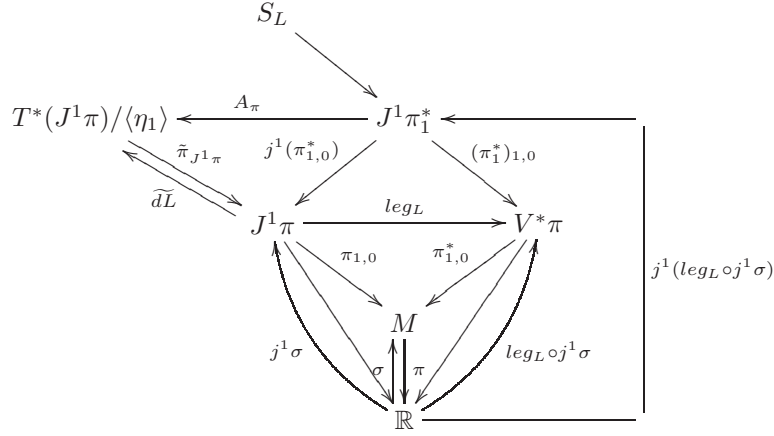
\begin{figure}[h]
\[
\xymatrix{
& S_{L}\ar[dr]&&\\
T^{*}(J^{1}\pi)/ \langle\eta_{1}\rangle\ar[rd]^{\tilde{\pi}_{J^{1}\pi}}&&J^{1}\pi^{*}_{1}\ar[ll]_{A_{\pi}}\ar[ld]_{j^{1}(\pi_{1,0}^{*})}\ar[rd]^{(\pi_{1}^{*})_{1,0}}&&\ar[ll]\\
& J^{1}\pi\ar@<1ex>[ul]^{\widetilde{dL}}\ar[rr]^{leg_{L}}\ar[rd]^{\pi_{1,0}}\ar[rdd]&& V^{*}\pi\ar[ld]_{\pi_{1,0}^{*}}\ar[ldd]&&\\
&&M\ar[d]^{\pi}&&&\\
&&\mathbb{R}\ar@{-}[rr]\ar[u]<1ex>^{\sigma}\ar@/^1pc/[uul]^{j^{1}\sigma}\ar@/_1pc/[uur]_{leg_{L}\circ j^{1}\sigma}&&\ar@{-}[uuu]_{j^{1}(leg_{L}\circ j^{1}\sigma)}&
 }
\]
\caption{\it The Lagrangian formalism in the restricted
Tulczyjew's triple}\label{fig2}
\end{figure}

Figure \ref{fig2}  illustrates the above situation

\subsection{The Hamiltonian formalism}
$\mbox{}$

Let $\mu: T^{*}M   \longrightarrow V^{*}\pi$ be the AV-bundle
associated with the fibration $\pi: M  \longrightarrow
\mathbb{R}$. $\mu$ defines a principal $ \mathbb{R}$-bundle.

We will denote by $V_{\mu}$ the infinitesimal generator of the action of $ \mathbb{R}$ on $T^{*}M$ and by

$$\emph{b}_{\Omega_{T^{*}M}}: T(T^{*}M) \longrightarrow  T^{*}(T^{*}M)$$
the vector bundle isomorphism (over the identity of $T^{*}M$)
induced by the canonical symplectic structure $\Omega_{T^{*}M}$ of
$T^{*}M$.

If $(t, q^{i}, p, p_{i}; \dot{t}, \dot{q}^{i}, \dot{p},
\dot{p}_{i})$ (respectively, $(t, q^{i}, p, p_{i}; p_{t},
p_{q^{i}}, p_{p}, p_{p_{i}})$ are local coordinates on $T(T^{*}M)$
(respectively, $T^{*}(T^{*}M)$), we have that
$$\emph{b}_{\Omega_{T^{*}M}}(t, q^{i}, p,  p_{i}; \dot{t}, \dot{q}^{i}, \dot{p}, \dot{p}_{i})= (t, q^{i}, p, p_{i}, - \dot{p}, -\dot{p}_{i}, \dot{t}, \dot{q}^{i}).$$
Now, if $\widehat{V_{\mu}}: T^{*}(T^{*}M) \longrightarrow
\mathbb{R}$ is the linear function on $T^{*}(T^{*}M)$ induced by
the vector field $V_{\mu}$, we can consider the affine subbundle
$\widehat{V_{\mu}}^{-1}(1)$ of $T^{*}(T^{*}M)$, that is,
$$\widehat{V_{\mu}}^{-1}(1) = \{\gamma \in T^{*}(T^{*}M) / \gamma(V_{\mu}(\pi_{T^{*}M}(\gamma)))=1\} $$
\noindent and the map $\varphi: \widehat{V_{\mu}}^{-1}(1)
\longrightarrow T(V^{*}\pi)$ defined by $\varphi = T\mu \circ
\emph{b}^{-1}_{\Omega_{T^{*}M}}$.

Since $\displaystyle{V_{\mu} =  \frac{\partial }{\partial p}}$ it
follows that $(t, q^{i}, p, p_{i}; p_{t}, p_{q^{i}}, p_{p_{i}})$
are local coordinates on $\widehat{V_{\mu}}^{-1}(1) $ and,
moreover,
$$\varphi(t, q^{i}, p, p_{i}; p_{t}, p_{q^{i}}, p_{p_{i}}) = (t, q^{i}, p_{i}; 1, p_{p_{i}}, -p_{q^{i}}).$$
Thus, $\varphi$ takes values in $J^{1}\pi_{1}^{*}$ and we can
consider the map
$$\varphi: \widehat{V_{\mu}}^{-1}(1) \longrightarrow J^{1}\pi_{1}^{*}.$$
The local expression of this map is
$$\varphi(t, q^{i}, p, p_{i}; p_{t}, p_{q^{i}}, p_{p_{i}}) =  (t, q^{i}, p_{i}; p_{p_{i}}, -p_{q^{i}}).$$
Therefore, it is clear that $\varphi$ is not a diffeomorphism. In
order to obtain a diffeomorphism, we will proceed as follows.

\emph{First Step}: The cotangent lift of the action of
$\mathbb{R}$ on $T^{*}M$ defines an action of $\mathbb{R}$ on
$T^{*}(T^{*}M)$. In fact, we have that
$$p'\cdot(t, q^{i}, p, p_{i}; p_{t}, p_{q^{i}},p_{p},  p_{p_{i}}) =(t, q^{i}, p+p', p_{i}; p_{t}, p_{q^{i}},p_{p},  p_{p_{i}})$$
for $p' \in \mathbb{R}$ and $(t, q^{i}, p, p_{i}; p_{t},
p_{q^{i}},p_{p},  p_{p_{i}}) \in T^{*}(T^{*}M)$.

It is obvious that the affine bundle $\widehat{V_{\mu}}^{-1}(1)$
is invariant under this action. Consequently, the space of orbits
of this action $\frac{\widehat{V_{\mu}}^{-1}(1)}{\mathbb{R}}$ is
an affine bundle over $V^{*}\pi$ which is modelled over the vector
bundle $\frac{\widehat{V_{\mu}}^{-1}(0)}{\mathbb{R}}$.
\begin{remark}{\rm The affine bundle $\frac{\widehat{V_{\mu}}^{-1}(1)}{\mathbb{R}}$ over $V^{*}\pi$
is identified with the \emph{phase bundle} $P\mu$ associated with
the AV-bundle $\mu: T^{*}M \longrightarrow V^{*}\pi$. The phase
bundle associated with an AV-bundle was introduced in
\cite{GrGrUr}.}\hfill{$\Diamond$}
\end{remark}
 Note that $\widehat{V_{\mu}}^{-1}(0)$  is just the annihilator of
the vertical bundle to $\mu: T^{*}M \longrightarrow V^{*}\pi$ and
that the quotient vector bundle
$\frac{\widehat{V_{\mu}}^{-1}(0)}{\mathbb{R}}$  is isomorphic to
$T^{*}(V^{*}\pi)$. So, the affine bundle $P\mu =
\frac{\widehat{V_{\mu}}^{-1}(1)}{\mathbb{R}}$  is modelled over
the vector bundle $T^{*}(V^{*}\pi)$.

Local coordinates on $P\mu =
\frac{\widehat{V_{\mu}}^{-1}(1)}{\mathbb{R}}$  are $(t, q^{i},
p_{i}; p_{t}, p_{q^{i}}, p_{p_{i}})$.

Moreover, there exists a smooth map $\overline{\varphi}: P\mu
\longrightarrow  J^{1}\pi_{1}^{*}$ such that the following diagram
\[
\xymatrix{
\widehat{V_{\mu}}^{-1}(1)\ar[rr]^{\varphi}\ar[d]_{\pi_{\widehat{V_{\mu}}^{-1}(1)}}&&J^{1}\pi^{*}_{1}\\
P\mu\ar[urr]_{\overline{\varphi}}&&
 }
\]
is commutative, where $\pi_{\widehat{V_{\mu}}^{-1}(1)}:
\widehat{V_{\mu}}^{-1}(1) \longrightarrow P\mu$ is the canonical
projection. The local expression of $\overline{\varphi}$ is
$$\overline{\varphi}(t, q^{i},  p_{i}; p_{t}, p_{q^{i}}, p_{p_{i}}) = (t, q^{i},  p_{i}; p_{p_{i}},- p_{q^{i}}).$$
Therefore, $\overline{\varphi}$ is a surjective submersion.

\emph{Second  Step}: Let $\pi_{1}^{*}: V^{*}\pi \longrightarrow
\mathbb{R}$ be the canonical projection. Then, the differential of
$\pi_{1}^{*}$ is a section of the vector bundle $\pi_{V^{*}\pi}:
T^{*}(V^{*}\pi) \longrightarrow V^{*}\pi$. Therefore, since $P\mu$
is an affine bundle modelled over $T^{*}(V^{*}\pi)$, we may
consider the quotient affine bundle $P\mu /\langle
d\pi^{*}_{1}\rangle$  over $V^{*}\pi$. $P\mu /\langle
d\pi^{*}_{1}\rangle$  is modelled over the quotient vector bundle
$T^{*}(V^{*}\pi) /\langle d\pi^{*}_{1}\rangle$ .

Local coordinates on $P\mu /\langle d\pi^{*}_{1}\rangle$  are $(t,
q^{i},  p_{i}; , p_{q^{i}}, p_{p_{i}})$.

Furthemore, there exits a smooth map $\tilde{\varphi}:P\mu
/\langle d\pi^{*}_{1}\rangle \longrightarrow J^{1}\pi_{1}^{*}$
such that the following diagram
\[
\xymatrix{
P\mu\ar[rr]^{\overline{\varphi}}\ar[d]_{\pi_{P\mu}}&&J^{1}\pi^{*}_{1}\\
P\mu /\langle d\pi^{*}_{1}\rangle\ar[urr]_{\tilde{\varphi}}&&
 }
\]
is commutative, where $\pi_{P\mu}: P\mu \longrightarrow P\mu
/\langle d\pi^{*}_{1}\rangle$ is the canonical projection. The
local expression of $\tilde{\varphi}$ is
$$\tilde{\varphi}(t, q^{i},  p_{i}; p_{q^{i}}, p_{p_{i}}) =(t, q^{i},  p_{i};  p_{p_{i}}, - p_{q^{i}}).$$
Consequently, $\tilde{\varphi}$ is a diffeomorphism.

We will denote by $\emph{b}_{\pi}: J^{1}\pi_{1}^{*}
\longrightarrow P\mu /\langle d\pi^{*}_{1}\rangle$ the inverse map
of $\tilde{\varphi}$, that is, $\emph{b}_{\pi}=
\tilde{\varphi}^{-1}$. Then, we have that
\begin{equation}\label{ec8}
\emph{b}_{\pi}(t, q^{i}, p_{i}, \dot{q}^{i}, \dot{p}_{i})= (t, q^{i}, p_{i},- \dot{p}_{i}, \dot{q}^{i}).
\end{equation}
Note that $\emph{b}_{\pi}$ is an affine bundle isomorphism over
the identity of $V^{*}\pi$.

The following diagram illustrates the situation
\[
\xymatrix{
J^{1}\pi^{*}_{1}\ar[rr]^{b_{\pi}}\ar[rd]_{(\pi_{1}^{*})_{1,0}}&&P\mu /\langle d\pi^{*}_{1}\rangle\ar[ld]^{\tilde{\pi}_{P\mu}}\\
& V^{*}\pi&&
 }
\]
Here $\tilde{\pi}_{P\mu}$ is the affine bundle projection.

$P\mu$ admits a canonical symplectic form $\Omega_{P\mu}$ (see
\cite{GrGrUr}). In fact, the local expression of $\Omega_{P\mu}$
is
$$\Omega_{P\mu}= dt\land dp_{t}+ dq^{i}\land dp_{q^{i}} + dp_{i} \land dp_{p_{i}}.$$
Let $\Lambda_{P\mu}$ be the Poisson structure on $P\mu$ associated
with $\Omega_{P\mu}$. Then,
$$\Lambda_{P\mu}=  \frac{\partial }{\partial t}\land \frac{\partial }{\partial p_{t}}+  \frac{\partial }{\partial q^{i}}\land \frac{\partial }{\partial p_{q^{i}}}+  \frac{\partial }{\partial p_{i}}\land \frac{\partial }{\partial p_{p_{i}}}.$$
On the other hand, the vertical lift $(d\pi_{1}^{*})^{v}$ to
$P\mu$ of the 1-form $d\pi_{1}^{*}$ on $V^{*}\pi$ generates the
vertical bundle to the canonical projection from $P\mu$ on $P\mu
/\langle d\pi^{*}_{1}\rangle$. Note that,
$$\displaystyle{(d\pi_{1}^{*})^{v} = \frac{\partial }{\partial p_{t}}}.$$

Thus, $ \displaystyle{\mathcal{L}_{(d\pi_{1}^{*})^{v} }
\Lambda_{P\mu}=0}$ and, therefore, $\Lambda_{P\mu}$ is projectable
to a Poisson structure $\widetilde{\Lambda}_{P\mu}$ on $P\mu
/\langle d\pi^{*}_{1}\rangle$.

The local expression of $\widetilde{\Lambda}_{P\mu}$ is
\begin{equation}\label{ec9}
\widetilde{\Lambda}_{P\mu}=  \frac{\partial }{\partial q^{i}}\land \frac{\partial }{\partial p_{q^{i}}}+  \frac{\partial }{\partial p_{i}}\land \frac{\partial }{\partial p_{p_{i}}}.
\end{equation}
Consequently, using  (\ref{ec7}),  (\ref{ec8}) and  (\ref{ec9}),
we prove the following result
\begin{theorem}\label{t1}
$\emph{b}_{\pi}$ is  anti-Poisson isomorphism between the Poisson manifolds $(J^{1}\pi_{1}^{*}, \Lambda_{J^{1}\pi_{1}^{*}})$ and  \\$( P\mu /\langle d\pi^{*}_{1}\rangle, \widetilde{\Lambda}_{P\mu})$.
\end{theorem}
Now, let $h: V^{*}\pi  \longrightarrow T^{*}M$  be a Hamiltonian
section and $F_{h}$ be the corresponding real function on $T^{*}M$
such that $V_{\mu}(F_{h})=1$. Then, one may define  a section of
the affine bundle $\widehat{V_{\mu}}^{-1}(1) \longrightarrow
T^{*}M$ as follows
$$\alpha \in T^{*}M  \longrightarrow dF_{h}(\alpha) \in \widehat{V_{\mu}}^{-1}(1).$$
This section is $\mathbb{R}$-equivariant. So, it induces a section
$dh: V^{*}\pi \longrightarrow P\mu$ of the phase bundle $P\mu$.
 We will denote by $\widetilde{dh}: V^{*}\pi  \longrightarrow  P\mu /\langle d\pi^{*}_{1}\rangle$
 the corresponding section of the affine bundle $ P\mu /\langle d\pi^{*}_{1}\rangle  \longrightarrow V^{*}\pi$.
 If the local expression of h  is
$$h(t, q^{i}, p_{i}) = (t, q^{i}, -H(t, q, p), p_{i}),$$
we have that
\begin{equation}\label{ec10}
\widetilde{dh}(t, q^{i}, p_{i}) = (t, q^{i}, p_{i} ;  \frac{\partial H }{\partial q^{i}} ,  \frac{\partial H}{\partial p_{i}}).
\end{equation}
Thus,
$$ (\tilde{\Lambda}_{P\mu}^{\sharp})^{-1}(T(\widetilde{dh}(V^{*}\pi))) = \Big{\langle} \Big{ \{ }dt, dp_{q^{j}} - \frac{\partial^{2}H}{\partial q^{i}\partial q^{j}}dq^{i} - \frac{\partial^{2}H}{\partial p_{i}\partial q^{j}} dp_{i}, \quad dp_{p_{j}} - \frac{\partial^{2}H}{\partial q^{i}\partial p_{j}}dq^{i} - \frac{\partial^{2}H}{\partial p_{i}\partial p_{j}}dp_{i} \Big{\} }\Big{\rangle},$$
\begin{eqnarray*}
 (\tilde{\Lambda}_{P\mu}^{\sharp}) \Big{(} T^{*}_{\widetilde{dh}(\alpha)}\Big{(} \frac{P\mu}{\langle d\pi_{1}^{*} \rangle} \Big{)}\Big{)} \cap T_{\widetilde{dh}(\alpha)}(\widetilde{dh}(V^{*}\pi)) &=& \Big{\langle} \Big{\{} \Big{(} \frac{\partial}{\partial q^{j}} + \frac{\partial^{2}H}{\partial q^{i}\partial q^{j}}\frac{\partial}{\partial p_{q^{i}}} + \frac{\partial^{2}H}{\partial q^{j}\partial p_{i}} \frac{\partial}{\partial p_{p_{i}}} \Big{)}_{\mid \widetilde{dh}(\alpha)},\\
  & & \quad \Big{(} \frac{\partial}{\partial p_{j}} + \frac{\partial^{2}H}{\partial q^{i}\partial p_{j}}\frac{\partial}{\partial p_{q^{i}}} + \frac{\partial^{2}H}{\partial p_{i}\partial p_{j}}\frac{\partial}{\partial p_{p_{i}}} \Big{)}_{\mid \widetilde{dh}(\alpha)} \Big{\}}\Big{\rangle},
\end{eqnarray*}
$\forall \alpha \in V^{*}\pi$.

Therefore,
$$\tilde{\Lambda}_{P\mu}(\alpha, \beta) = 0, \quad \forall \alpha, \beta \in (\tilde{\Lambda}_{P\mu}^{\sharp})^{-1}(T(\widetilde{dh}(V^{*}\pi))),$$
$$dim \Big{ (}T_{\widetilde{dh}(\alpha)}(\widetilde{dh}(V^{*}\pi)) \cap (\tilde{\Lambda}_{P\mu}^{\sharp})\Big{(}T_{\widetilde{dh}(\alpha)}\Big{(} \frac{P\mu}{\langle d\pi_{1}^{*} \rangle} \Big{)}\Big{)} \Big{)}= \frac{ dim \Big{(}\tilde{\Lambda}_{P\mu}^{\sharp}\Big{(}T^{*}_{\widetilde{dh}(\alpha)} \Big{(} \frac{P\mu}{\langle d\pi_{1}^{*} \rangle } \Big{)} \Big{)}\Big{)} }{2} = 2n,$$
$\forall \alpha \in V^{*}\pi.$

This implies that $\widetilde{dh}(V^{*}\pi)$ is a Lagrangian
submanifold of the Poisson manifold $\Big{(} \displaystyle
\frac{P\mu}{\langle d\pi_{1}^{*} \rangle}, \tilde{\Lambda}_{P\mu}
\Big{)}.$

So, from Theorem \ref{t1}, it follows that $S_{h} =
\emph{b}_{\pi}^{-1}(\widetilde{dh}(V^{*}\pi))$ is also a
Lagrangian submanifold of the Poisson manifold $(J^{1}\pi_{1}^{*},
\Lambda_{J^{1}\pi_{1}^{*}})$.

 On the other hand, if $R_{h}$ is
the  Reeb vector field of the cosymplectic structure $(\Omega_{h},
\eta_{1}^{*})$ on $V^{*}\pi$ (see subsection \ref{Ham-formalism})
then, using  (\ref{Rh}), (\ref{ec8}) and (\ref{ec10}), we deduce
that
$$S_{h}= R_{h}(V^{*}\pi).$$
Consequently, since the integral curves of $R_{h}$ are the
solutions of the Hamilton equations for the Hamiltonian section
$h$, we obtain the following result.
\begin{theorem}
\begin{enumerate}
\item Let $\tau: \mathbb{R} \longrightarrow V^{*}\pi$ be a local
section of the fibration $\pi_{1}^{*}: V^{*}\pi  \longrightarrow  \mathbb{R}$.
Then,  $\tau$ is a solution of the Hamilton equations for h if and only
if
$$\emph{b}_{\pi}^{-1} \circ \widetilde{dh} \circ \tau = j^{1}\tau.$$
\item The local equations which define to $S_{h}$ as a Lagrangian submanifold of
the Poisson ma\-ni\-fold $(J^{1}\pi_{1}^{*},
\Lambda_{J^{1}\pi_{1}^{*}})$ are just the Hamilton equations for
$h$.
\end{enumerate}
\end{theorem}

\begin{figure}[h]
\[
\xymatrix{
&S_{h}\ar[dl]&\\
J^{1}\pi_{1}^{*}\ar[rr]^{\emph{b}_{\pi}}\ar[dr]^{(\pi_{1}^{*})_{1,0}}&&P\mu/\langle d\pi_{1}^{*}\rangle\ar[dl]_{\tilde{\pi}_{P\mu}}\\
&V^{*}\pi\ar@<-1ex>[ur]_{\widetilde{dh}}\ar[d]^{\pi_{1,0}^{*}}&\\
&M\ar[ld]^{\pi}&\\
\mathbb{R}\ar[uuu]^{j^{1}\tau}\ar[uur]^{\tau}&&
}
\]
\caption{\it The Hamiltonian formalism in the restricted
Tulczyjew's triple}\label{fig3}
\end{figure}
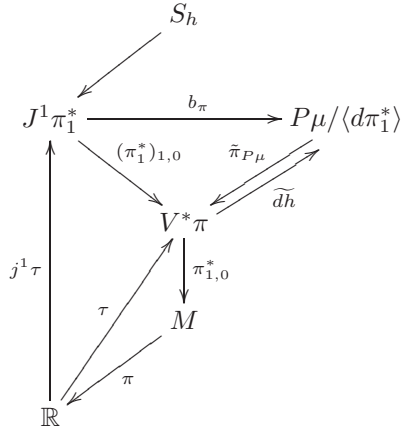

Figure \ref{fig3} illustrates the situation

\subsection{The equivalence between the Lagrangian and Hamiltonian formalism}
$\mbox{}$

Let $L: J^{1}\pi  \longrightarrow  \mathbb{R}$ be an hyperregular
Lagrangian function. Then, the restricted Legendre transformation
$leg_{L}: J^{1}\pi \longrightarrow V^{*}\pi$ is a global
diffeomorphism and we may consider the Euler-Lagrange vector field
$R_{L}$ on $J^{1}\pi$. Note that, since $
leg_{L}^{*}(\eta_{1}^{*}) = \eta_{1}$ and $\eta_{1}(R_{L}) = 1$,
it follows that $Tleg_{L}(R_{L}(J^{1}\pi)) \subseteq
J^{1}\pi_{1}^{*}.$

Moreover, using (\ref{ec3}), (\ref{R-L}), (\ref{ec5}) and
(\ref{dl}), we deduce
\begin{lemma}\label{l1}
The following relation holds
$$A_{\pi} \circ Tleg_{L} \circ R_{L} = \widetilde{dL}.$$
\end{lemma}

\begin{figure}[h]
\[
\xymatrix{
& S_{L}\ar[dr]&=&S_{h}\ar[dl]\\
T^{*}(J^{1}\pi)/ \langle\eta_{1}\rangle\ar[rd]^{\tilde{\pi}_{J^{1}\pi}}&&J^{1}\pi^{*}_{1}\ar[ll]_{A_{\pi}}\ar[rr]^{\emph{b}_{\pi}}\ar[ld]_{j^{1}(\pi_{1,0}^{*})}\ar[rd]^{(\pi_{1}^{*})_{1,0}}&&P\mu/\langle d\pi^{*}_{1} \rangle\ar[dl]_{\tilde{\pi}_{P\mu}}&\\
& J^{1}\pi\ar@<1ex>[ul]^{\widetilde{dL}}\ar[rr]^{leg_{L}}&& V^{*}\pi\ar@<-1ex>[ur]_{\tilde{dh}}&&
 }
\]
\caption{\it The restricted Tulczyjew's triple for time-dependent
Mechanics}\label{fig4}
\end{figure}
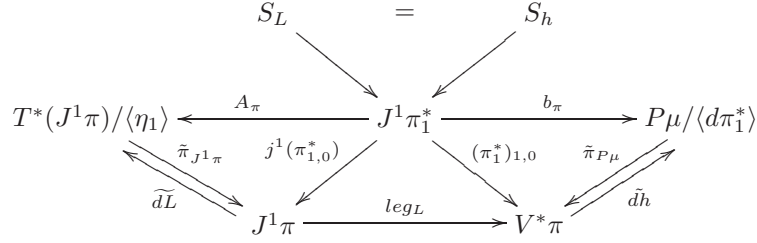

Now, denote by $h: V^{*}\pi  \longrightarrow  T^{*}M$ the
Hamiltonian section associated with the hyperregular Lagrangian
function L, that is,
$$h= Leg_{L} \circ leg_{L}^{-1}, $$
$Leg_{L}: J^{1}\pi  \longrightarrow  T^{*}M$ being the extended
Legendre transformation. Then, using Lemma \ref{l1} and since
$Tleg_{L} \circ R_{L} = R_{h} \circ leg_{L}$, we prove the
following result.
\begin{theorem}
The Lagrangian submanifolds $S_{L} =
A_{\pi}^{-1}(\widetilde{dL}(J^{1}\pi))$ and
$S_{h}=R_{h}(V^{*}\pi)$ of the Poisson manifold
$(J^{1}\pi_{1}^{*}, \Lambda_{J^{1}\pi_{1}^{*}})$ are equal.
\end{theorem}
The previous result may be considered as the expression of the
equivalence between the Lagrangian formalism and the restricted
Hamiltonian formalism in the Lagrangian submanifold setting.
Figure \ref{fig4} illustrates the situation

\section{Extended Tulczyjew's triple}

\subsection{The Lagrangian formalism}
\setcounter{equation}{0}

$\mbox{}$

Let $\tilde{\pi}_{M} : T^{*}M \longrightarrow \mathbb{R}$ be the
fibration from $T^{*}M$ on $\mathbb{R}$. We consider the space
$J^{1}\tilde{\pi}_{M}$ of 1-jets of local sections of
$\tilde{\pi}_{M} : T^{*}M \longrightarrow \mathbb{R}$. As we know,
there exists a natural embedding from $ J^{1}\tilde{\pi}_{M}$ in
$T(T^{*}M)$, which we will denote by $ j: J^{1}\tilde{\pi}_{M}
\longrightarrow T(T^{*}M)$.

On the other hand, we can consider the 1-jet prolongation
$j^{1}\pi_{M}: J^{1}\tilde{\pi}_{M} \longrightarrow J^{1}\pi$ of
the bundle map $\pi_{M}: T^{*}M  \longrightarrow M$.

Then, we may define a smooth map
$$\widetilde{A_{\pi}} : J^{1}\tilde{\pi}_{M} \longrightarrow T^{*}(J^{1}\pi)$$
as follows:

Let $\tilde{z}$ be a point of $J^{1}\tilde{\pi}_{M}$ and $A_{M}:
T(T^{*}M)  \longrightarrow T^{*}(TM)$ be the canonical Tulczyjew
diffeomorphism. Then, $A_{M}(j(\tilde{z})) \in T_{v}^{*}(TM)$,
with $v \in J^{1}\pi$. Indeed, if $(t, q^{i}, p, p_{i})$ are local
coordinates on $T^{*}M$, we have that  $(t, q^{i}, p, p_{i};
\dot{q}^{i}, \dot{p}, \dot{p}_{i})$ are local coordinates on
$J^{1}\tilde{\pi}_{M}$ and

$$A_{M}(j(\tilde{z})) = (t, q^{i}, 1, \dot{q}^{i}; \dot{p}, \dot{p}_{i}, p, p_{i}).$$
 Thus, $A_{M}(j(\tilde{z})) \in T_{v}^{*}(TM)$, with $v \in
J^{1}\pi$. In fact, $v= (j^{1}\pi_{M})(\tilde{z})$.

Now, we define
$$ \displaystyle{\widetilde{A_{\pi}}(\tilde{z}) = A_{M}(j(\tilde{z}))_{\mid T_{j^{1}\pi_{M}(\tilde{z})}(J^{1}\pi)} \in T^{*}_{(j^{1}\pi_{M}(\tilde{z})}(J^{1}\pi).}$$
Therefore, it follows that
\begin{equation}\label{atilde}
\widetilde{A_{\pi}}(t, q^{i}, p, p_{i}; \dot{q}^{i}, \dot{p}, \dot{p}_{i}) = (t, q^{i}, \dot{q}^{i}; \dot{p}, \dot{p}_{i}, p_{i}).
\end{equation}

Consequently, $\widetilde{A_{\pi}}$ is a surjective submersion.
$\widetilde{A_{\pi}}$ is called \emph{the canonical Tulczyjew
fibration a\-sso\-cia\-ted with $\pi$}.
\begin{remark}\label{fibers-orbits}{\rm $\widetilde{A_{\pi}}$ is
the bundle projection of a principal $\mathbb{R}$-bundle. In fact,
if we consider the tangent lift of the principal action of
$\mathbb{R}$ on $T^*M$, we have an action of $\mathbb{R}$ on
$T(T^*M)$. The local expression of this action is
\[
p'\cdot (t, q^i, p, p_i; \dot{t}, \dot{q}^i, \dot{p}, \dot{p}_i) =
(t, q^i, p+p', p_i; \dot{t}, \dot{q}^i, \dot{p}, \dot{p}_i)
\]
for $p' \in \mathbb{R}$ and $(t, q^i, p, p_i; \dot{t}, \dot{q}^i,
\dot{p}, \dot{p}_i) \in T(T^*M)$.

Thus, it is clear that the submanifold $J^1\tilde{\pi}_{M}$ of
$T(T^*M)$ is invariant under the previous action and, from
(\ref{atilde}), it follows that the fibers of
$\widetilde{A_{\pi}}$ are just the orbits of the action of
$\mathbb{R}$ on $J^1\tilde{\pi}_{M}$. }\hfill{$\Diamond$}
\end{remark}
Next, we will denote by $\Omega_{M}$ the canonical symplectic
structure of $T^{*}M$ and by $\Omega_{M}^{c}$ the complete lift of
$\Omega_{M}$ to $T(T^{*}M)$. $\Omega_{M}^{c}$ defines a symplectic
structure on $T(T^{*}M)$ and $j^{*}(\Omega_{M}^{c}) =
\Omega_{J^{1}\tilde{\pi}_{M}}$ is a presymplectic form on
$J^{1}\tilde{\pi}_{M}$.

In fact,  the local expressions of these forms are
$$\Omega_{M}^{c} = dt\land d\dot{p } +  d\dot{t}\land dp+ dq^{i}\land d\dot{p}_{i}+  d\dot{q}^{i}\land dp_{i},$$
and
\begin{equation}\label{pre}
 \Omega_{J^{1}\tilde{\pi}_{M}} = dt\land d\dot{p } + dq^{i}\land d\dot{p}_{i}+  d\dot{q}^{i}\land dp_{i}.
\end{equation}

Thus, $  \Omega_{J^{1}\tilde{\pi}_{M}}$ is a presymplectic form of
corank 1 and the kernel of $  \Omega_{J^{1}\tilde{\pi}_{M}}$ is
generated by the restriction to $J^{1}\tilde{\pi}_{M}$ of the
complete lift ${(V_{\mu})}^{c}$ of $V_{\mu}$ to $T(T^{*}M)$. Note
that,
\begin{equation}\label{nucleo}
{(V_{\mu})}^{c} = \frac{\partial}{\partial p} \quad and  \quad
ker(T\widetilde{A_{\pi}}) = \langle \{ {(V_{\mu})}^{c} \}\rangle.
\end{equation}

On the other hand, let $ \Omega_{J^{1}\pi}$ be the canonical
symplectic structure of $T^{*}(J^{1}\pi)$ . Then, if $(t, q^{i},
\dot{q}^{i}; p_{t}, \linebreak p_{q^{i}}, p_{\dot{q}^{i}})$ are
local coordinates on $T^{*}(J^{1}\pi)$, we have that
\begin{equation}\label{sym}
\Omega_{J^{1}\pi}= dt\wedge dp_{t} + dq^{i}\wedge dp_{q^{i}} + d\dot{q}^{i}\wedge dp_{\dot{q}^{i}}.
\end{equation}
Therefore, using (\ref{atilde}), (\ref{pre}) and (\ref{sym}), we
deduce the following result.
\begin{theorem}\label{tT}
The canonical Tulczyjew fibration associated with $\pi$ is a
presymplectic map between the presymplectic manifolds
$(J^{1}\tilde{\pi}_{M},  \Omega_{J^{1}\tilde{\pi}_{M}})$ and
$(T^{*}(J^{1}\pi),  \Omega_{J^{1}\pi})$, that is,
$$\widetilde{A_{\pi}}^{*}(\Omega_{J^{1}\pi}) = \Omega_{J^{1}\tilde{\pi}_{M}}.$$
\end{theorem}
Now, let $ L: J^{1}\pi \longrightarrow  \mathbb{R}$ be a
Lagrangian function. Then, it is well-known that $dL(J^{1}\pi)$ is
a Lagrangian submanifold of the symplectic manifold
$(T^{*}(J^{1}\pi),  \Omega_{J^{1}\pi})$. Consequently, using
(\ref{nucleo}) and  Theorem \ref{tT}, we obtain that
$\widetilde{S_{L}}= \widetilde{A_{\pi}}^{-1}(dL(J^{1}\pi))$ also
is a Lagrangian submanifold of the presymplectic  manifold
$(J^{1}\tilde{\pi}_{M},  \Omega_{J^{1}\tilde{\pi}_{M}})$.

Moreover, if $\sigma$ is a local section of $\pi: M
\longrightarrow  \mathbb{R}$ then, from (\ref{ec3}) and
(\ref{atilde}), we deduce that
$$\widetilde{A_{\pi}} \circ j^{1} (Leg_{L} \circ (j^{1}\sigma)(t)) =\Big{( }t, q^{i}(t), E_{L}(j^{1}\sigma(t)), \frac{\partial L}{\partial\dot{q}^{i}}((j^{1}\sigma)(t));  \frac{dq^{i}}{dt}, \frac{d(E_{L}\circ j^{1}\sigma)}{dt}, \frac{d}{dt}\Big{(}\frac{\partial L}{\partial \dot{q}^{i}} \circ j^{1}\sigma\Big{)}
\Big{)}$$
where $E_{L} = L - \dot{q}^{i}\frac{\partial L}{\partial
\dot{q}^{i}}$ and $Leg_{L}: J^{1}\pi \longrightarrow T^{*}M$ is
the extended Legendre transformation (see (\ref{ec3})).

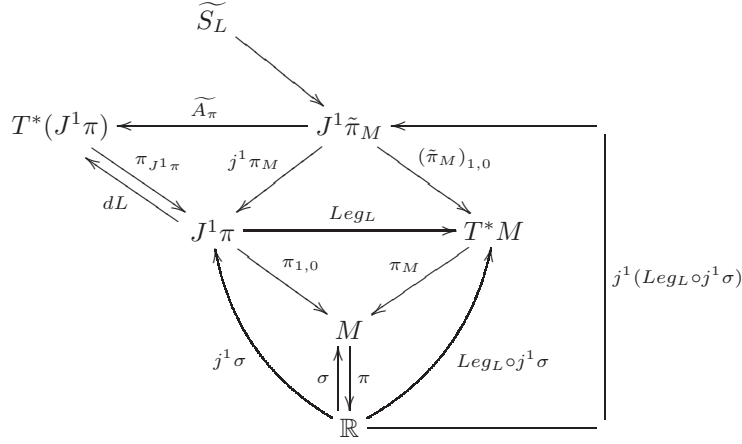
\begin{figure}[h]
\[
\xymatrix{
&\widetilde{ S_{L}}\ar[dr]&&\\
T^{*}(J^{1}\pi)\ar[rd]^{\pi_{J^{1}\pi}}&&J^{1}\tilde{\pi}_{M}\ar[ll]_{\widetilde{A_{\pi}}}\ar[ld]_{j^{1}\pi_{M}}\ar[rd]^{{(\tilde{\pi}_{M})}_{1,0}}&&\ar[ll]\\
& J^{1}\pi\ar@<1ex>[ul]^{dL}\ar[rr]^{Leg_{L}}\ar[rd]^{\pi_{1,0}}&& T^{*}M\ar[ld]_{\pi_{M}}&&\\
&&M\ar[d]^{\pi}&&&\\
&&\mathbb{R}\ar@{-}[rr]\ar[u]<1ex>^{\sigma}\ar@/^1pc/[uul]^{j^{1}\sigma}\ar@/_1pc/[uur]_{Leg_{L}\circ j^{1}\sigma}&&\ar@{-}[uuu]_{j^{1}(Leg_{L}\circ j^{1}\sigma)}&
 }
\]
\caption{\it The Lagrangian formalism in the extended Tulczyjew's
triple}\label{fig5}
\end{figure}

We remark that for a solution $\sigma$ of the Euler-Lagrange equations for L, we  have that
$$ \frac{d(E_{L}\circ j^{1}\sigma)}{dt} = \frac{\partial L}{\partial t} \circ j^{1}\sigma.$$
Using the above facts, one may prove the following result.
\begin{theorem}
\begin{enumerate}
\item A section $\sigma:  \mathbb{R} \longrightarrow M $ is a solution of the Euler-Lagrange equations
for L if and only if
$$dL \circ j^{1}\sigma = \widetilde{A_{\pi}} \circ j^{1}(Leg_{L} \circ
j^{1}\sigma).$$
\item The local equations which define to $\widetilde{S_{L}}$ as a Lagrangian submanifold
of the presymplectic ma\-ni\-fold $(J^{1}\tilde{\pi}_{M},
\Omega_{J^{1}\tilde{\pi}_{M}})$ are just the Euler-Lagrange
equations for L.
\end{enumerate}
\end{theorem}
Figure \ref{fig5} illustrates the situation


\subsection{The Hamiltonian formalism}
$\mbox{}$

Let $\tilde{\pi}_{M}: T^{*}M \longrightarrow \mathbb{R}$ be the
fibration from $T^{*}M$ on $\mathbb{R}$. Recall that
$J^{1}\tilde{\pi}_{M}$ is the space of 1-jets of local sections of
$\tilde{\pi}_{M} : T^{*}M \longrightarrow \mathbb{R}$ and that j
is the  natural embedding from $ J^{1}\tilde{\pi}_{M}$ in
$T(T^{*}M)$.

Then, we may define a  map
$$\widetilde{\emph{b}_{\pi}}: J^{1}\tilde{\pi}_{M} \longrightarrow T^{*}(T^{*}M) $$
as follows:

Let $\tilde{z}$ be a point of $J^{1}\tilde{\pi}_{M}$ and
$\emph{b}_{M}: T(T^{*}M) \longrightarrow T^{*}(T^{*}M)$ the vector
bundle isomorphism (over the identity of $T^{*}M$) induced by the
canonical symplectic structure $\Omega_{M}$ of $T^{*}M$. Then,
$\widetilde{\emph{b}_{\pi}}(\tilde{z}) =
\emph{b}_{M}(j(\tilde{z})) \in T^{*}_{\alpha}(T^{*}M)$, with
$\alpha \in T^{*}M$. In fact, if $(t, q^{i}, p, p_{i})$ are local
coordinates on $T^{*}M$, we have that  $(t, q^{i}, p, p_{i};
\dot{q}^{i}, \dot{p}, \dot{p}_{i})$ are local coordinates on
$J^{1}\tilde{\pi}_{M}$ and
$$\widetilde{\emph{b}_{\pi}}(t, q^{i}, p, p_{i}; \dot{q}^{i}, \dot{p}, \dot{p}_{i}) =  (t, q^{i}, p, p_{i}; -\dot{p}, - \dot{p}_{i}, 1,
\dot{q}^{i}).$$
From the last equation, we observe that the map
$\widetilde{\emph{b}_{\pi}}$ takes values on the affine subbundle
$\widehat{V_{\mu}}^{-1}(1)$ of $T^{*}(T^{*}M)$. For this reason,
we can consider the map
$$\widetilde{\emph{b}_{\pi}}: J^{1}\tilde{\pi}_{M} \longrightarrow \widehat{V_{\mu}}^{-1}(1)$$
\noindent which in local coordinates is given by
\begin{equation}\label{btilde}
\widetilde{\emph{b}_{\pi}}(t, q^{i}, p, p_{i}; \dot{q}^{i}, \dot{p}, \dot{p}_{i}) = (t, q^{i}, p, p_{i}; -\dot{p}, - \dot{p}_{i}, \dot{q}^{i}).
\end{equation}
Consequently, $\widetilde{\emph{b}_{\pi}}$ is a diffeomorphism.
\begin{remark}\label{equivariant}{\rm
If we consider the cotangent lift of the principal action of
$\mathbb{R}$ on $T^*M$, we have an action of $\mathbb{R}$ on
$T^*(T^*M)$. The local expression of this action is
\[
p'\cdot (t, q^i, p, p_i; p_t, p_{q^i}, p_p, p_{p_{i}}) = (t, q^i,
p+p', p_i; p_t, p_{q^i}, p_p, p_{p_{i}})
\]
for $p'\in \mathbb{R}$ and $(t, q^i, p, p_i; p_t, p_{q^i}, p_p,
p_{p_{i}}) \in T^*(T^*M)$.

Thus, it is clear that the affine subbundle
$\widehat{V_{\mu}}^{-1}(1)$ of $T^*(T^*M)$ is invariant under this
action. Moreover, if we consider the natural action of
$\mathbb{R}$ on $J^{1}\tilde{\pi}_{M}$ (see Remark
\ref{fibers-orbits}) then, from (\ref{btilde}), it follows that
the diffeomorphism $\widetilde{\emph{b}_{\pi}}$ is equivariant.
}\hfill{$\Diamond$}
\end{remark}
Next, we will denote by $\Omega_{T^{*}M}$ the canonical symplectic
structure on $T^{*}(T^{*}M)$ and by $\Phi_{
\widehat{V_{\mu}}^{-1}(1)}$ the 2-form on $
\widehat{V_{\mu}}^{-1}(1)$ defined by
$$\Phi_{ \widehat{V_{\mu}}^{-1}(1)} = i_{ \widehat{V_{\mu}}^{-1}(1)}^{*}(\Omega_{T^{*}M}),$$
where $i_{ \widehat{V_{\mu}}^{-1}(1)}:  \widehat{V_{\mu}}^{-1}(1)
\longrightarrow T^{*}(T^{*}M)$ is the canonical inclusion.

The local expressions of these forms are
$$\Omega_{T^{*}M} = dt \land dp_{t} + dq^{i} \land dp_{q^{i}} + dp \land dp_{p} + dp_{i}\land dp_{p_{i}}, $$
and
\begin{equation}\label{omegaA}
\Phi_{ \widehat{V_{\mu}}^{-1}(1)} = dt \land dp_{t} + dq^{i} \land
dp_{q^{i}} + dp_{i}\land dp_{p_{i}}.
\end{equation}
Thus,  $\Phi_{ \widehat{V_{\mu}}^{-1}(1)}$ is a presymplectic form
of corank 1  and the kernel of $\Phi_{ \widehat{V_{\mu}}^{-1}(1)}$
is generated by the restriction to $\widehat{V_{\mu}}^{-1}(1)$ of
the complete lift $(V_{\mu})^{*c}$ of $V_{\mu}$ to
$T^{*}(T^{*}M)$. Note that $(V_{\mu})^{*c}$ is the Hamiltonian
vector field of the linear function $\widehat{V_{\mu}}: T^*(T^*M)
\to \mathbb{R}$ and, therefore,
\begin{equation}\label{vu}
(V_{\mu})^{*c}= \frac{\partial}{\partial p}.
\end{equation}
Consequently, using (\ref{pre}), (\ref{btilde}) and
(\ref{omegaA}), we deduce the following result.
\begin{theorem}\label{thepre}
$\widetilde{\emph{b}_{\pi}}:J^{1}\tilde{\pi}_{M} \longrightarrow
\widehat{V_{\mu}}^{-1}(1)$ is an anti-presymplectic isomorphism
between the presymplectic manifolds $(J^{1}\tilde{\pi}_{M},
\Omega_{J^{1}\tilde{\pi}_{M}})$ and $( \widehat{V_{\mu}}^{-1}(1),
\Phi_{ \widehat{V_{\mu}}^{-1}(1)})$, that is,
$$\widetilde{\emph{b}_{\pi}}^{*}(\Phi_{ \widehat{V_{\mu}}^{-1}(1)}) = -  \Omega_{J^{1}\tilde{\pi}_{M}}.$$
\end{theorem}
Now, let $h: V^{*}\pi  \longrightarrow T^{*}M$ be a Hamiltonian
section and $F_{h}: T^{*}M  \longrightarrow  \mathbb{R}$ be the
correspon\-ding real $ C^{\infty}$-function on $T^{*}M$ satisfying
$V_{\mu}(F_{h})=1$ (see section 3.2). Then, it is clear that
$dF_{h}(T^{*}M) \subseteq \widehat{V_{\mu}}^{-1}(1) \subseteq
T^{*}(T^{*}M)$.

Denote by $i_{dF_{h}(T^{*}M)}: dF_{h}(T^{*}M) \longrightarrow
\widehat{V_{\mu}}^{-1}(1)$ the canonical inclusion.

Since $dF_{h}(T^{*}M)$ is a Lagrangian submanifold of
$T^{*}(T^{*}M)$ and $\Phi_{ \widehat{V_{\mu}}^{-1}(1)} = i_{
\widehat{V_{\mu}}^{-1}(1)}^{*}(\Omega_{T^{*}M})$, we deduce that
\begin{equation}\label{n1}
i_{dF_{h}(T^{*}M)}^{*}(\Phi_{ \widehat{V_{\mu}}^{-1}(1)}) = 0.
\end{equation}
On the other hand, using (\ref{vu}), it is easy to prove that the
restriction of $(V_{\mu})^{*c}$ to $dF_{h}(T^{*}M)$ is tangent to
$dF_{h}(T^{*}M)$. Thus,
\begin{equation}\label{n2}
Ker\Big{(}\Phi_{ \widehat{V_{\mu}}^{-1}(1)}(dF_{h}(\alpha))\Big{)}
\subseteq T_{dF_{h}(\alpha)}(dF_{h}(T^{*}M)), \quad \forall \alpha
\in T^{*}M.
\end{equation}
Therefore, from (\ref{n1}) and (\ref{n2}), we obtain that
$dF_{h}(T^{*}M)$ is a Lagrangian submanifold of the presymplectic
manifold $( \widehat{V_{\mu}}^{-1}(1),  \Phi_{
\widehat{V_{\mu}}^{-1}(1)})$ (see Definition \ref{lspm}).

Consequently, using Theorem \ref{thepre}, it follows that
$\widetilde{S_{h}} =
\widetilde{\emph{b}_{\pi}}^{-1}(dF_{h}(T^{*}M))$ is also a
Lagrangian submanifold of the presymplectic manifold
$(J^{1}\tilde{\pi}_{M},  \Omega_{J^{1}\tilde{\pi}_{M}})$.

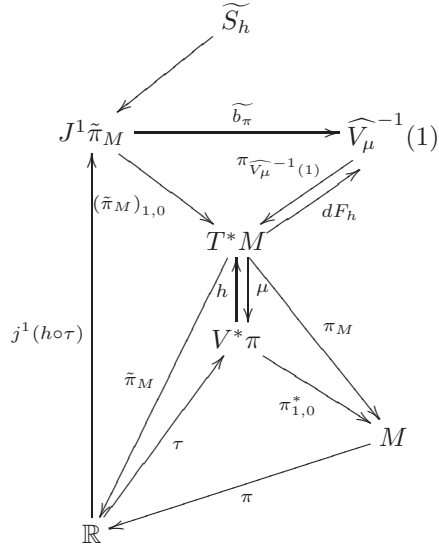
\begin{figure}[h]
\[
\xymatrix{
&\widetilde{S_{h}}\ar[dl]&\\
J^{1}\tilde{\pi}_{M}\ar[dr]_{{(\tilde{\pi}_{M})}_{1,0}}\ar[rr]^{\widetilde{\emph{b}_{\pi}}}&&\widehat{V_{\mu}}^{-1}(1)\ar[dl]_{\pi_{\widehat{V_{\mu}}^{-1}(1)}}\\
&T^{*}M\ar[dddl]_{\tilde{\pi}_{M}}\ar@<1ex>[d]^{\mu}\ar@<-1ex>[ur]_{dF_{h}}\ar[ddr]^{\pi_{M}}&\\
&V^{*}\pi\ar[u]^{h}\ar[dr]_{\pi_{1,0}^{*}}&\\
&&M\ar[dll]^{\pi}\\
\mathbb{R}\ar[uur]_{\tau}\ar[uuuu]^{j^{1}(h\circ \tau)}&&
}
\]
\caption{\it The Hamiltonian formalism in the extended Tulczyjew's
triple}\label{fig6}
\end{figure}

Next, suppose that  $\tau: \mathbb{R}\longrightarrow V^{*}\pi$ is
a section of $\pi^{*}_{1}: V^{*}\pi \longrightarrow \mathbb{R}$.
Then, we have that
$$(dF_{h} \circ h \circ \tau)( \mathbb{R})  \subseteq  \widehat{V_{\mu}}^{-1}(1)$$
(see (\ref{Fh})). Moreover, if $\tau$ is a solution of the
Hamilton equations then, from (\ref{a2}), we deduce that
$$\frac{d(H\circ \tau)}{dt} = \frac{\partial H}{\partial t} \circ \tau.$$
Using these facts and (\ref{btilde}), we may prove the following
result.
\begin{theorem}
\begin{enumerate}
\item  A section  $\tau: \mathbb{R}\longrightarrow V^{*}\pi$ is a solution of Hamilton
equations for h if and only if
 $$\widetilde{\emph{b}_{\pi}} \circ j^{1}(h \circ \tau) = dF_{h} \circ h \circ \tau.$$
\item The local equations which define to $\widetilde{S_{h}}$ as a Lagrangian submanifold
of  $J^{1}\tilde{\pi}_{M}$ are just the Hamilton equations for h.
\end{enumerate}
\end{theorem}
Figure \ref{fig6} illustrates the situation

\subsection{The equivalence between the Lagrangian and Hamiltonian formalism}
$\mbox{}$

Let $L: J^{1}\pi  \longrightarrow  \mathbb{R}$ be an hyperregular
Lagrangian function. Then, the restricted Legendre transformation
$leg_{L}: J^{1}\pi \longrightarrow V^{*}\pi$ is a global
diffeomorphism and we may consider the Euler-Lagrange vector field
$R_{L}$ on $J^{1}\pi$.

Moreover, using (\ref{ec3}), (\ref{R-L}) and (\ref{atilde}), we
deduce
\begin{lemma}\label{l2}
The following relation holds
$$\widetilde{A_{\pi}} \circ TLeg_{L} \circ R_{L} = dL,$$
where $Leg_{L}: J^1\pi \to T^*M$ is the extended Legendre
transformation.
\end{lemma}
Now, denote by $h: V^{*}\pi  \longrightarrow  T^{*}M$ the
Hamiltonian section associated with the hyperregular Lagrangian
function L, that is,
$$h= Leg_{L} \circ leg_{L}^{-1}.$$
\begin{theorem}
The Lagrangian submanifolds $\widetilde{S_{L}}=
\widetilde{A_{\pi}}^{-1}(dL(J^{1}\pi)) $ and $\widetilde{S_{h}} =
\widetilde{\emph{b}_{\pi}}^{-1}(dF_{h}(T^{*}M))$ of the
presymplectic  manifold $(J^{1}\tilde{\pi}_{M},
\Omega_{J^{1}\tilde{\pi}_{M}})$ are equal.
\end{theorem}
\underline{Proof:} Let $\tilde{z}$ be a point of
$\widetilde{S_{L}}$. Then, since $\pi_{J^{1}\pi} \circ
\widetilde{A_{\pi}} = j^{1}\pi_{M}$, it follows that
$$\widetilde{A_{\pi}}(\tilde{z}) = dL((j^{1}\pi_{M})(\tilde{z})).$$
Thus, using Lemma \ref{l2} and the fact that $R_{L}$ and
$\mathcal{H}^{\Omega_{M}}_{F_{h}}$ are $Leg_{L}$-related, we
deduce that
$$\widetilde{A_{\pi}}(\tilde{z}) = \widetilde{A_{\pi}}(\mathcal{H}^{\Omega_{M}}_{F_{h}}(Leg_{L}(j^{1}\pi_{M})(\tilde{z}))) =
\widetilde{A_{\pi}}(\widetilde{\emph{b}_{\pi}}^{-1}(dF_{h}(Leg_{L}(j^{1}\pi_{M})(\tilde{z})))).$$
Therefore, from Remark \ref{fibers-orbits}, we obtain that there
exists a unique $p\in \mathbb{R}$ such that
\[
\widetilde{\emph{b}_{\pi}}(p\cdot \tilde{z}) =
dF_{h}(Leg_{L}((j^1\pi_{M})(\tilde{z}))).
\]
Here, $\cdot$ denotes the action of $\mathbb{R}$ on
$J^1\tilde{\pi}_{M}$.

Consequently, using Remarks \ref{princ-connec} and
\ref{equivariant}, it follows that
\[
\widetilde{\emph{b}_{\pi}}(\tilde{z}) = dF_{h}((-p) \cdot
Leg_{L}((j^1\pi_{M})(\tilde{z}))) \in dF_{h}(T^*M).
\]
So, $\tilde{z} \in \widetilde{\emph{b}_{\pi}}^{-1}(dF_{h}(T^*M)) =
\widetilde{S_h}$. This implies that $\widetilde{S_{L}}\subseteq
\widetilde{S_{h}}$.

Proceeding in a similar way, one may prove that $\widetilde{S_{h}}
\subseteq \widetilde{S_{L}}$. \hfill{$\square$}

The previous result may be considered as the expression of the
equivalence between the Lagrangian and extended  Hamiltonian
formalism in the Lagrangian submanifold setting.

\begin{figure}[h]
\[
\xymatrix{
& \widetilde{S_{L}}\ar[dr]&=&\widetilde{S_{h}}\ar[dl]\\
T^{*}(J^{1}\pi)\ar[rd]^{\pi_{J^{1}\pi}}&&J^{1}\tilde{\pi}_{M}\ar[ll]_{\widetilde{A_{\pi}}}\ar[rr]^{\widetilde{\emph{b}_{\pi}}}\ar[ld]_{j^{1}\pi_{M}}\ar[rd]^{{(\tilde{\pi}_{M})}_{1,0}}&&\widehat{V_{\mu}}^{-1}(1)\ar[dl]_{\pi_{\widehat{V_{\mu}}^{-1}(1)}}&\\
& J^{1}\pi\ar@<1ex>[ul]^{dL}\ar[rr]^{Leg_{L}}&& T^{*}M\ar@<-1ex>[ur]_{dF_{h}}&&
 }
\]
\caption{\it The extended Tulczyjew's triple for time-dependent
Mechanics}\label{fig7}
\end{figure}
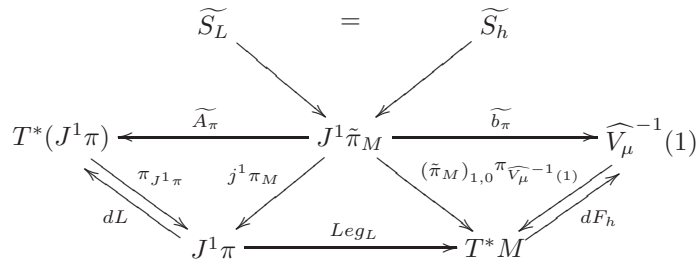

Figure \ref{fig7} illustrates the situation

\vspace{-10pt}

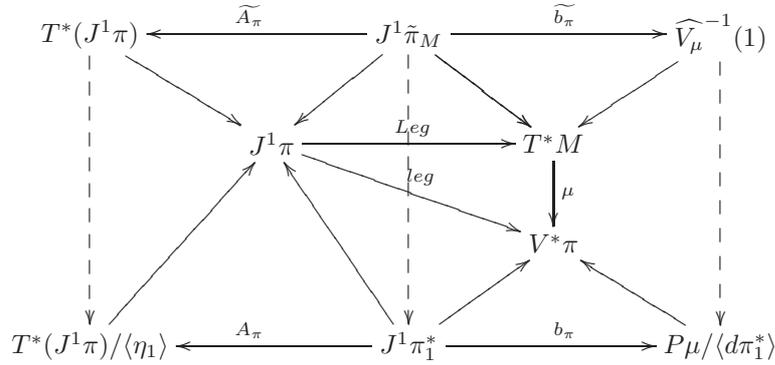
\begin{figure}[h]
\[
\xymatrix{
T^{*}(J^{1}\pi)\ar@{-->}[ddd]\ar[rd]&&J^{1}\tilde{\pi}_{M}\ar[rr]^{\widetilde{b_{\pi}}}\ar@{-->}[ddd]\ar[rd]\ar[ll]_{\widetilde{A_{\pi}}}\ar[ld]\ar[rd]&&\widehat{V_{\mu}}^{-1}(1)\ar[ld]\ar@{-->}[ddd]\\
& J^{1}\pi\ar[rr]^{Leg}\ar[rrd]^{leg}&& T^{*}M\ar[d]^{\mu}&&\\
&&&V^{*}\pi&&\\
T^{*}(J^{1}\pi)/ \langle\eta_{1}\rangle\ar[ruu]&&J^{1}\pi^{*}_{1}\ar[rr]^{b_{\pi}}\ar[ll]_{A_{\pi}}\ar[luu]\ar[ru]&&P\mu /\langle d\pi^{*}_{1}\rangle\ar[lu]\\
 }
\]
\caption{\it The restricted and  extended Tulczyjew's triples for
time-dependent Mechanics}\label{fig8}
\end{figure}

Finally, Figure \ref{fig8} describes both triples. The extended
Tulczyjew triple is on the top of the diagram and the restricted
Tulczyjew triple is on the bottom.

\vspace{-10pt}

\section{Conclusions and future work}

Using the geometry of presymplectic and Poisson manifolds a new Tulczyjew
triple for time-dependent Mechanics is discussed. More precisely,
we present two Tulczyjew triples. The first one is adapted to the
restricted Hamiltonian formalism for time-dependent mechanical
systems and the second one is adapted to the extended Hamiltonian
formalism. Our construction solves some problems and deficiences
of previous approaches.

It would be interesting to extend the ideas and results contained
in this paper for classical field theories of first order. For
this purpose,  a suitable higher order generalization of a  presymplectic (Poisson)
structure must be used. This will be the subject of a
forthcoming paper.

Other Tulczyjew triples for classical field theories of  first
order have been proposed by several authors (see \cite{Gr,
LeMaSa}).




\begin{thebibliography}{99}
  \let\\, \newcommand{\by}[1]{\textsc{#1}\\}
  \newcommand{\title}[1]{\textsl{#1}\\}
  \newcommand{\vol}[1]{\textbf{#1}}
  \newcommand{\info}[1]{\textrm{#1}.}

\bibitem{AbMa} \by{R Abraham and JE Marsden}
\title{Foundations of Mechanics} \info{Second edition, revised and
enlarged., Benjamin/Cummings, Reading, Mass., 1978}

\bibitem{BuRa} \by{H Burzstyn and O Radko} \title{Gauge equivalence
of Dirac structures and symplectic groupoids} \info{Ann. Inst.
Fourier (Grenoble) \vol{53} (2003), 309--337}

\bibitem{Co} \by{TJ Courant} \title{Dirac manifolds} \info{Trans. Amer. Math.
Soc. \vol{319} (1990), 631--661}

\bibitem{Gr} \by{K Grabowska} \title{Lagrangian and Hamiltonian
formalism in Field Theory: a simple model} \info{Preprint,
arXiv:1005.2753}

\bibitem{GrGrUr} \by{K Grabowska, J Grabowski and P Urba\'nski}
\title{AV-differential Geometry: Poisson and Jacobi structures}
\info{Journal of Geometry and Physics \vol{52} (2004), 398--446}

\bibitem{GrGrUr1} \by{K Grabowska, J Grabowski and P Urba\'nski}
\title{AV-differential geometry: Euler-Lagrange equations}
\info{J. Geom. Phys. \vol{57} (2007), 1984--1998}

\bibitem{GrUr} \by{J Grabowski and P Urba\'nski}
\title{Tangent lifts of Poisson and related structures}
\info{J. Phys. A: Math. Gen. \vol{28} (1995), 6743--6777}

\bibitem{IgMaPaSo} \by{D Iglesias, JC Marrero, E Padr\'on and D Sosa}
\title{Lagrangian submanifolds and dynamics on Lie affgebroids} \info{Rep.
Math. Phys. \vol{38} (2006), 385--436}

\bibitem{LeLa} \by{M de Le\'on and EA Lacomba}
\title{Lagrangian submanifolds and higher-order mechanical systems} \info{J. Phys.
A: Math. Gen. \vol{22} (1989), 3809--3820}

\bibitem{LeMaMa} \by{M de Le\'on, J Mar{\'\i}n-Solano and JC Marrero}
\title{The constraint algorithm in the jet formalism} \info{
Differential Geom. Appl. \vol{6} (3) (1996), 275--300}

\bibitem{LeMa} \by{M de Le\'on and JC Marrero}
\title{Constrained time-dependent Lagrangian systems and Lagrangian
submanifolds} \info{J. Math. Phys. \vol{34} (1993), 622--644}

\bibitem{LeMaSa} \by{M de Le\'on, D Mart{\'\i}n de Diego and
A Santamar{\'\i}a-Merino} \title{Tulczyjew triples and lagrangian
submanifolds in classical field theories} \info{in Applied
Differential Geometry and Mechanics, Editors W Sarlet and F
Cantrijn, Univ. of Gent, Gent, Academia Press, 2003, 2147}

\bibitem{LeRo} \by{M de Le\'on and PR Rodrigues}
  \title{Methods of Differential Geometry in Analytical Mechanics}
  \info{North Holland Math.\ Series \vol{152} (Amsterdam, 1996)}

\bibitem{LiMa} \by{P Libermann and ChM Marle} \title{Symplectic Geometry and
Analytical Mechanics} \info{Kluwer, Dordrecht, (1987)}

\bibitem{Sa} \by{DJ Saunders} \title{The geometry of jet bundles}
\info{London Math. Soc., Lecture Note Series, \vol{142} Cambridge
Univ. Press, (1989)}

\bibitem{Tu1} \by{W Tulczyjew} \title{Les sous-vari\'{e}t\'{e}s lagrangiennes et la
dynamique hamiltonienne} \info{C.R. Acad. Sci. Paris \vol{283}
(1976), 15--18}

\bibitem{Tu2} \by{W Tulczyjew} \title{Les sous-vari\'{e}t\'{e}s lagrangiennes et la
dynamique lagrangienne} \info{C.R. Acad. Sci. Paris \vol{283}
(1976), 675--678}

\bibitem{Uc} \by{K Uchino} \title{Lagrangian calculus on Dirac manifolds}
\info{J. Math. Soc. Japan \vol{57} (3) (2005), 803--825}


\end{thebibliography}
\end{document}